\documentclass[preprint]{imsart}
\usepackage{xr}
\RequirePackage[OT1]{fontenc}
\RequirePackage{amsthm,amsmath}
\RequirePackage{natbib,color}
\RequirePackage[colorlinks,citecolor=blue,urlcolor=blue]{hyperref}

\startlocaldefs
\numberwithin{equation}{section}
\theoremstyle{plain}

\endlocaldefs

\usepackage{psfrag,
graphicx, 
epsf,yhmath,
euscript,amsfonts,amsmath,latexsym,amssymb,mathrsfs
}
\usepackage[active]{srcltx}
\usepackage[a4paper]{geometry}

\newtheorem{lemma}{Lemma}
\newtheorem{theorem}{Theorem}

\newtheorem{example}{Example}
\newtheorem{assumption}{Assumption}
\newtheorem{remark}{Remark}
\newtheorem{definition}{Definition}

\newcommand{\cC}{{\cal C}}

\newcommand{\cF}{{\cal F}}

\newcommand{\cL}{{\cal L}}

\newcommand{\cR}{{\cal R}}

\newcommand{\bC}{\mathbb C}

\newcommand{\bE}{\mathbb E}

\newcommand{\bL}{{\mathbb L}}

\newcommand{\bN}{{\mathbb N}}

\newcommand{\bP}{{\mathbb P}}

\newcommand{\bR}{{\mathbb R}}

\newcommand{\bZ}{{\mathbb Z}}
\newcommand{\sF}{{\mathscr F}}
\newcommand{\sG}{{\mathscr G}}
\newcommand{\sH}{{\mathscr H}}

\newcommand{\sL}{{\mathscr L}}
\newcommand{\sM}{{\mathscr M}}

\newcommand{\sS}{{\mathscr S}}

\newcommand{\rd}{\mathrm{d}}

\renewcommand{\kappa}{\varkappa}

\newcommand{\epr}{\hfill\hbox{\hskip 4pt
                \vrule width 5pt height 6pt depth 1.5pt}\vspace{0.5cm}\par}
\begin{document}

\begin{frontmatter}
\title{Density deconvolution under general assumptions 
on the distribution 
of  
measurement errors
}
\runtitle{Density deconvolution  under general assumptions}

\begin{aug}
\author{\fnms{Denis} \snm{Belomestny}
\ead[label=e1]{denis.belomestny@uni-due.de}} \and
\author{\fnms{Alexander} \snm{Goldenshluger}
\ead[label=e2]{goldensh@stat.haifa.ac.il}}

\runauthor{D. Belomestny and A. Goldenshluger}

\affiliation{Duisburg-Essen University\thanksmark{m1}, 
Higher School of Economics\thanksmark{m3} 
\\
and University of Haifa\thanksmark{m2}}

\address{Faculty of Mathematics\\ 
Duisburg-Essen University
\\
Thea-Leymann-Str. 9
\\
D-45127 Essen
\\
Germany \\
\printead{e1}\\
\phantom{E-mail:\ }
}

\address{Department of Statistics
\\ University of Haifa
\\ Haifa  3498838 
\\
Israel
\\
\printead{e2}}

\address{National University 
\\
Higher School of Economics
\\
11 Pokrovsky Bulvar, Pokrovka Complex
\\
Moscow, Russia
}

\end{aug}

\begin{abstract} 
In  this paper we study   the problem 
of density deconvolution under general assumptions on the measurement error distribution. Typically 
deconvolution estimators are constructed using Fourier transform techniques, and 
it is assumed that 
the  characteristic function of 
the measurement errors does not have zeros 
on the real line. This assumption is rather strong and is not fulfilled 
in many cases of interest.  In this paper we develop a 
methodology for constructing optimal density deconvolution estimators in the general setting that covers 
vanishing  and non--vanishing  characteristic functions of the measurement errors.
We derive upper bounds on the risk of the proposed estimators and 
provide  sufficient conditions under which zeros of the corresponding characteristic function have no effect on estimation accuracy.
Moreover, we show that the  derived conditions are also necessary in some
specific problem instances.
\end{abstract}
 

\end{frontmatter}

\section{Introduction}\label{sec:introduction}

\subsection{Problem formulation and background}
The problem of density deconvolution is formulated as follows.
Suppose that
we observe a random sample $Y_1,\ldots,Y_n$ generated from the model
\begin{equation*}
 Y_j=X_j+\epsilon_j,\;\;\;j=1,\ldots, n,
\end{equation*}
where $X_1,\ldots, X_n$ are i.i.d. random variables with density $f$, and the measurement errors
$\epsilon_1,\ldots,\epsilon_n$ are i.i.d. random variables with known distribution $G$. Furthermore, assume that 
$\epsilon_1,\ldots, \epsilon_n$ are independent of $X_1,\ldots, X_n$.  
Then the probability density $f_Y$ of $Y=X+\epsilon$ is given by convolution
\begin{equation}\label{eq:convolution}
 f_Y(y)= [f\ast \rd G](y)=\int_{-\infty}^\infty f(y-x) \rd G(x).
\end{equation}
The goal is to estimate $f$ from the observations $Y_1,\ldots, Y_n$. 
\par 
An estimator 
$\tilde{f}$ of $f$ is 
%
a measurable 
function of observations $(Y_1, \ldots, Y_n)$, and the accuracy of $\tilde{f}$ 
is measured
by the maximal  risk
\[
 \cR_{n, \Delta} [\tilde{f}; \sF]= \sup_{f\in \sF} 
 \Big\{\bE_f \Delta^2(\tilde{f},f) \Big\}^{1/2},
 %
\]
where 
$\Delta(\cdot, \cdot)$ is a loss function, 
$\sF$ is  a class of density functions,  and  $\bE_f$ 
is
the expectation with respect to the 
probability measure $\bP_f$ of observations $Y_1,\ldots, Y_n$ when the density of $X$ is $f$.
In this paper we will be interested in estimating $f$ at a single point $x_0\in \bR$ and in the $\bL_2$--norm; this
corresponds to the loss 
functions $\Delta_{x_0}(f_1,f_2):=|f_1(x_0)-f_2(x_0)|$  and
$\Delta_2(f_1, f_2):=\|f_1-f_2\|_2=\big\{\int_{-\infty}^\infty |f_1(x)-f_2(x)|^2\rd x\big\}^{1/2}$,
respectively.
The minimax risk is then defined by 
\[
\cR^*_{n, \Delta}[\sF]=\inf_{\tilde{f}}\cR_{n, \Delta}[\tilde{f}; \sF],
\]
where $\inf$ is taken over all possible estimators. An estimator
$\tilde{f}_*$ is called {\em rate--optimal} 
if $\cR_{n, \Delta}[\tilde{f}_*;\sF]\asymp O(\cR_{n, \Delta}^*[\sF])$ as $n\to\infty$, and our goal is  to construct rate--optimal estimators
for natural functional classes of densities.
\par 
The problem of density deconvolution has been extensively studied in the literature; see, e.g.,    
\cite{Carroll-Hall88}, \cite{Zhang90}, \cite{Fan91}, \cite{Butucea-T1}, \cite{Lounici-Nickl}, \cite{Comte-Lacour13} and 
\cite{Lepski-Willer19}.
We also refer to
the book of \cite{Meister}, where many additional references can be found. 
\par 
Deconvolution estimators are usually constructed using Fourier transform techniques, and 
the majority of results in the existing literature assumes that 
the characteristic function of the measurement errors has no zeros on the real line.
Specifically, let $\widehat{g}$ denote the bilateral Laplace transform of $G$, 
\[ 
\widehat{g}(z):= \int_{-\infty}^\infty e^{-z x}\rd G(x),
\] 
with $\widehat{g}(i\omega)$ being 
the characteristic function of the measurement errors.
The standard assumptions on $\widehat{g}$  in the density deconvolution problem  
 are the following:
\begin{itemize}
\item[(A)] $\widehat{g}(i\omega)$ does not vanish, that is, $|\widehat{g}(i\omega)|\ne 0$  for all  $\omega \in \bR$;
\item[(B)] $|\widehat{g}(i\omega)|$ decreases in an appropriate way as $|\omega|\to \infty$: for some $\gamma>0$
\begin{itemize}
\item[(B1)] 
$|\widehat{g}(i\omega)|\asymp |\omega|^{-\gamma}$ as $|\omega|\to\infty$, 
\item [or]
\item[(B2)] 
$|\widehat{g}(i\omega)|\asymp \exp\{-c|\omega|^\gamma\}$ as $|\omega|\to\infty$ with $c>0$.
\end{itemize}
\end{itemize} 
\par 
The setting under  conditions (A)--(B1)
is usually referred to as   the {\em case of  smooth} measurement error densities, while conditions (A)--(B2) correspond to the so--called 
{\em super--smooth case}.
Under assumption~(A) the achievable estimation accuracy is determined 
by  the rate at which $|\widehat{g}(i\omega)|$ decreases as $|\omega|\to\infty$, and by 
smoothness of the density $f$. In particular, it is well known that in {\em the smooth case} 
for the H\"older class $\sH_\alpha(A)$ and for the Sobolev class $\sS_\alpha(A)$ of regularity $\alpha$ one has 
\begin{equation}
\label{eq:opt-int}
\cR^*_{n, \Delta_{x_0}}[\sH_\alpha(A)]\asymp n^{-\alpha/(2\alpha+2\gamma+1)}, \quad
\cR^*_{n, \Delta_2}[\sS_\alpha(A)]\asymp n^{-\alpha/(2\alpha+2\gamma+1)}
\end{equation}
as $n\to\infty$;
see, e.g., \cite{Zhang90} and \cite{Fan91}.
The definitions of classes $\sH_\alpha(A)$ and $\sS_\alpha(A)$ are presented in Section~\ref{sec:upper-bound}. 
In all what follows we will refer to the rate $n^{-\alpha/(2\alpha+2\gamma+1)}$ as the {\em standard rate of convergence}.
\par
It is worth noting that
the condition~(A) is  rather restrictive and excludes many settings of interest. 
This condition does not hold if distribution of 
the measurement errors is compactly supported. 
For instance, if $g$ is a uniform density on $[-1,1]$ then $\widehat{g}(i\omega)=\sin \omega/\omega$, and $\widehat{g}(i\omega)$
vanishes at $\omega=\pi k$, $k=\pm 1, \pm 2,\ldots$.
Another typical situation in which condition~(A) is violated   is the case of measurement errors having
discrete distributions.
In general, if $\widehat{g}(i\omega)$  has zeros,  
the standard Fourier--transform--based estimation methods are not directly applicable. 
This fact raises the  following natural questions. 
\begin{itemize}
\item[(i)]  How to construct 
rate-optimal estimators in the case 
when the assumption (A) does not hold, that is, $\widehat{g}(i\omega)$ has zeros, and  what is the best achievable rate of convergence 
under these  circumstances ?  
\item[(ii)] Under which conditions on \(f\)
one can achieve  the standard rates of convergence \eqref{eq:opt-int} without assuming (A) ?
\end{itemize}
\par 
The existing literature  contains only 
partial and 
fragmentary answers to  questions~(i) and~(ii).
\cite{Devroye89} constructed a consistent estimator of~$f$ under assumption that $|\widehat{g}(i\omega)|\ne 0$ for {\em almost all} $\omega$. The proposed estimator is a certain modification of the standard Fourier--transform--based kernel density estimator. 
\cite{Hall-etal01} consider the setting with 
the uniform measurement error density $g$ and develop an estimator under assumption 
that density $f$ is a compactly supported. Other works dealing with the uniform density 
deconvolution are \cite{GroenJong03}
and \cite{FeuerKim08}. The first cited paper
assumes that $X$ is  non--negative,
and shows that for a class of twice continuously differentiable densities, 
the pointwise risk  of the proposed estimators converges to 
zero at the 
standard rate corresponding to $\gamma=1$.
\cite{FeuerKim08} studied  the problem of estimating densities from 
Sobolev functional classes with the $\bL_2$--risk; they show that the standard rate of convergence with $\gamma=1$ can be 
achieved in this setting provided 
that $f$ has two bounded moments.  These results demonstrate that, in the problem with 
uniformly distributed measurement errors and
under 
the aforementioned assumptions on~$f$,  zeros of the 
characteristic function of~$\epsilon$ have no effect on the minimax rate of convergence. 
\par 
\cite{Hall-Meister} and 
\cite{Meister-1}
considered density deconvolution problem
with an oscillating Fourier transform $\widehat{g}(i\omega)$ that vanishes periodically. They proposed several 
modifications of the standard Fourier--transform--based estimators, considered the $\bL_2$--risk and showed that for  
certain nonparametric classes of probability densities,  zeros of the 
characteristic
function $\widehat{g}(i\omega)$ {\em do affect} the rate of convergence.  
\cite{Delaigle-Meister} demonstrated that if the density 
to be estimated has a finite left endpoint, then it can be estimated with the standard rate as in the case where 
$\widehat{g}(i\omega)$ does not have zeros.
\cite{Meister-Neumann}
considered a  setting where $\widehat{g}(i\omega)$ may have zeros, 
but there are two observations 
of the same variable $X$ with  independent measurement errors. In this setting   zeros of 
$\widehat{g}(i\omega)$  
have no influence on the rate of convergence.
\par 
The existing results in the literature leave open a fundamental question  about  
construction  of  rate--optimal density deconvolution estimators  under general assumptions on 
$\widehat{g}(i\omega).$
Specifically,
it is not clear whether  and under which conditions
zeros of $\widehat{g}(i\omega)$ have no influence on the minimax rates of convergence.
\par 
The current  paper addresses the aforementioned  issues. 
First, we develop a methodology for constructing optimal density deconvolution estimators
under general conditions on the measurement error distribution.  These conditions cover settings with vanishing and non--vanishing 
characteristic functions of the measurement errors, and the proposed methodology treats all these settings
in a unified way.  
The estimation methods we propose are based 
on  the 
Laplace transform. In this sense they  generalize  standard  
Fourier--transform--based estimation techniques  commonly used in the literature on density deconvolution.  
Second, we derive upper bounds on the risk of the proposed estimators and 
provide  sufficient conditions on $f$ under which the standard rate of convergence can be  achieved under general assumptions 
on~$\widehat{g}$.
In particular, we prove that if,  
in addition to  the smoothness restriction $f\in \sH_\alpha(A)$ or $f\in \sS_\alpha(A)$,  density $f$ has 
bounded moments of  a sufficiently large order, 
then the standard rate of convergence can be achieved even without Assumption~(A). 
The required number of bounded moments 
is 
characterized  in terms of  a sequence of coefficients 
({\em zero set sequence}) which, in turn, is determined  
by the geometry of zeros of $\widehat{g}(i\omega)$.
Third, 
we specialize our general methodology to specific problem instances in which the zero set sequences can be
explicitly calculated.  Last but not least, it is also shown 
that the derived sufficient moment conditions 
are also necessary  
in order to guarantee  the standard rate of convergence in absence of (A) for some specific problem instances.  
\par 
The rest of the paper is organized as follows.
In Section~\ref{sec:idea} we present a general idea for construction of proposed estimators. Section~\ref{sec:assumptions}
introduces assumptions on the distribution of the measurement errors and presents 
examples of distributions satisfying these assumptions. Section~\ref{sec:representation}
discusses construction of the estimator kernel and develops its infinite series representation.
In Section~\ref{sec:upper-bound} we define the estimator and present upper bounds on its risk.
Settings corresponding to 
specific problem instances are discussed in Section~\ref{sec:special}, and 
lower bounds showing necessity of moment conditions  are  presented in Section~\ref{sec:lower-bound}.
Some concluding remarks are brought in Section~\ref{sec:conclusion}. 
Proofs of all results are given in 
Appendix.
\subsection{Notation}
For a generic locally integrable function $\phi$ the bilateral Laplace transform is 
defined by
\[
 \widehat{\phi}(z)=\cL[\phi; z]:= \int_{-\infty}^\infty \phi(t) e^{-zt}\rd t.
\]
The Laplace transform $\widehat{\phi}(z)$ is  an analytic function in the 
convergence region $\Sigma_\phi$ of the above integral which, in general, is
a vertical strip:
\[
\Sigma_\phi:= \{z\in \bC:\sigma_\phi^-<{\rm Re}(z)<\sigma_\phi^+\} \;\;\;\hbox{for some}\;\;\; 
-\infty \leq \sigma_\phi^- < \sigma_{\phi}^+\leq \infty.
\]
The convergence region can degenerate to
a vertical line 
$\Sigma_\phi:=\{z\in \bC: {\rm Re}(z)=\sigma_\phi\}$, $\sigma_\phi\in\bR$, 
in the complex plane. If $\phi$ is a probability density, 
then 
the imaginary axis always belongs to $\Sigma_\phi$, that is,
$\{z:{\rm Re}(z)=0\}\subseteq \Sigma_\phi$, and 
\[
\widehat{\phi}(i\omega)=\cF[\phi; \omega]:=\int_{-\infty}^\infty \phi(t) e^{-i\omega t}\rd t, \;\;\omega\in \bR 
\]
is the characteristic function (the Fourier transform of $\phi$). This degenerate case  corresponds to distributions whose characteristic function cannot be analytically continued  to a strip  around the imaginary axes in the complex plane. 
The inverse Laplace transform is given by the formula
\[
 \phi(t)=\frac{1}{2\pi i}\int_{s-i\infty}^{s+i\infty} \widehat{\phi}(z) e^{zt}\rd z =
 \frac{1}{2\pi} \int_{-\infty}^\infty \widehat{\phi}(s+i\omega) e^{(s+i\omega) t} \rd \omega,\;\;\;
 s\in (\sigma_\phi^-, \sigma_\phi^+).
\]
The uniqueness property of the bilateral Laplace transform states that  if $\widehat{\phi}_1(z)=\widehat{\phi}_2(z)$
in a common strip of convergence ${\rm Re}(z)\in (\sigma_{\phi_1}^-, \sigma_{\phi_1}^+)\cap (\sigma_{\phi_2}^-, \sigma_{\phi_2}^+),$
then $\phi_1(t)$ is equal to $\phi_2(t)$ for almost all $t$ 
\cite[Theorem~6b]{Widder46}.
\section{General idea for estimator construction}\label{sec:idea}
Let $G$ be the measurement error distribution function
with the corresponding Laplace transform
\[ 
\widehat{g}(z)= \int_{-\infty}^\infty e^{-zt} \rd G(t)
\]
whose convergence region is denoted $\Sigma_g$.
%
%
Throughout the paper we suppose that
$\Sigma_g$ is a vertical strip in the complex plane, 
$\Sigma_g=\{z: \sigma_g^{-}<{\rm Re}(z)<\sigma_g^+\}$
with $ \sigma_g^-$ and $\sigma_g^+$ satisfying $\sigma_g^-<0<\sigma_g^+$ (see Assumption~\ref{as:general} in Section~\ref{sec:assumptions}). 
As it was discussed above, if $\widehat{g}(z)$ has zeros on the imaginary  axis in the complex plane, 
then the usual Fourier--transform--based methods are not directly applicable. We 
will be mainly interested in this case.
\subsection{Linear functional strategy}
The construction of our estimators follows the so-called {\em linear functional strategy} that is frequently used
for solving ill--posed inverse problems [see, e.g., \cite{goldberg1979amethod} and \cite{anderssen1980ontheuse}].
In the  context of the density deconvolution problem the main idea of the strategy is as follows. 
We want to find two kernels, say, $K$ and $L$ with the following properties:  
\begin{itemize}
 \item[(i)] integral $\int K(x-x_0) f(x)\rd x$ approximates ``well'' the value $f(x_0)$; 
 \item[(ii)] the kernel $L$  is  related to the kernel $K$ via the equation:
 \begin{equation}\label{eq:lin-funct}
 \int K(x) f(x)\rd x=\int L(y) f_Y(y)\rd y.
 \end{equation}
\end{itemize}
Under  conditions~(i) and~(ii), the obvious estimator of $f(x_0)$ from the observations 
$Y_1, \ldots, Y_n$ is the empirical estimator $\tilde{f}$ of the integral on the right hand side
of (\ref{eq:lin-funct}):
\[
\tilde{f}(x_0):=\frac{1}{n}\sum_{i=1}^n L(Y_i).
\] 
\par
Let $K:\bR\to \bR$ be a kernel with standard properties that will be specified later. For $h>0$ denote  $K_h(\cdot):=(1/h)K(\cdot/h)$. 
Assume that $K$ has bounded support so that 
$\widehat{K}(z)$ is an entire function, that is, $\Sigma_K=\bC$.
Furthermore, 
assume that 
there exist real numbers $\kappa_g^-$ and $\kappa_g^+$ satisfying $\sigma_g^-\leq \kappa_g^-<0<\kappa_g^+\leq \sigma_g^+$, 
such that 
\begin{equation}\label{eq:g-}
 \widehat{g}(-z)\ne 0,\;\;\;\;\forall z\in S_g:=\big\{z: {\rm Re}(z)\in (\kappa_g^-, 0)\cup (0, \kappa_g^+)\big\}.
\end{equation}
In words,  $S_g$ is the union of two open strips (with the imaginary axis as the boundary), where the function $\widehat{g}(-z)$ does not have zeros. 
Therefore we can define 
\begin{equation}\label{eq:L-hat}
 \widehat{L}_h (z):=\frac{\widehat{K}(zh)}{\widehat{g}(-z)},\;\;\;z\in S_g,
\end{equation}
and this function is 
analytic in $S_g$. 
Let
\begin{eqnarray}
 L_{s,h}(t)&:=& 
 \frac{1}{2\pi i}\int_{s-i\infty}^{s+i\infty} 
 \frac{\widehat{K}(zh)}{\widehat{g}(-z)} e^{z t} \rd z
 \nonumber
 \\
 &=&
 \frac{1}{2\pi} \int_{-\infty}^\infty \frac{\widehat{K}((s+i\omega)h)}{\widehat{g}(-s-i\omega)} e^{(s+i\omega)t}\rd \omega\;\;\;
 \label{eq:L}
\end{eqnarray}
with \(s\in (\kappa_g^-, 0)\cup (0, \kappa_g^+).\) Observe that the kernel $L_{s,h}$ is defined by the inverse Laplace transform of the function 
$\widehat{L}_h(z)=\widehat{K}(zh)/\widehat{g}(-z)$, and the 
denominator of the integrand in (\ref{eq:L}) does not vanish as 
$s\in (\kappa_g^-, 0)\cup (0, \kappa_g^+)$.
If the integral on the right hand side of (\ref{eq:L}) is absolutely convergent then (\ref{eq:L})
defines the same function $L_{s,h}$ for any value of $s\in (\kappa_g^-,0)$ or 
$s\in (0,\kappa_g^+)$. 
In other words, depending on the sign of $s$, the equation (\ref{eq:L})
defines two different functions which will be denoted by $L_{+, h}(t)$ and $L_{-,h}(t)$, respectively. 
The  estimator of $f(x_0)$ is then defined by 
\begin{equation}\label{eq:f-sh}
 \tilde{f}_{s,h}(x_0)= \frac{1}{n}\sum_{j=1}^n L_{s,h}(Y_j-x_0), \;\;\;s\in (\kappa_g^-, 0) \cup (0, \kappa_g^+).
\end{equation}
The parameters $s$ and $h$ 
will be specified in the sequel. 
\subsection{Relationship between kernels $K_h$ and $L_{s,h}$} 
The following lemma demonstrates that (\ref{eq:lin-funct}) holds for  
the kernels $K_h(\cdot):= (1/h)K(\cdot/h)$ and $L_{s,h}$ given by
(\ref{eq:L}).
\begin{lemma}\label{lem:K-L}
Suppose that for any  $s\in (\kappa_g^-, 0)\cup (0, \kappa_g^+)$
the integral on the right hand side of (\ref{eq:L}) is absolutely convergent, and 
\[
 \int_{-\infty}^\infty |L_{s, h}(y-x_0)| f_Y(y)\rd y <\infty;
\]
then for any $x_0$
\begin{equation}\label{eq:L-K}
 \int_{-\infty}^\infty L_{s, h}(y-x_0) f_Y(y)\rd y = \int_{-\infty}^\infty \frac{1}{h} K\bigg(\frac{x-x_0}{h} \bigg) f(x)\rd x.
\end{equation}
\end{lemma}

Note that relation (\ref{eq:L-K}) holds for both kernels 
$L_{+, h}$ and $L_{-, h}$ corresponding to $s\in (0, \kappa_g^+)$ and 
$s\in (\kappa_g^-, 0),$ respectively.  
Thus, both $L_{+,h}$ or $L_{-,h}$ can be used in the estimator construction.
\begin{remark}
A naive approach towards estimator construction 
could be based on a direct application of 
the Laplace transform inversion formula. 
In particular, 
(\ref{eq:convolution}) implies that $\widehat{f}_Y(z)=\widehat{f}(z)\widehat{g}(z)$. The empirical estimator of $\widehat{f}_Y(z)$ can be constructed
in the standard way using 
the available data $Y_1,\ldots, Y_n$; then a division by $\widehat{g}(z)$ with 
proper regularization and application of the inverse Laplace transform formula
 yields an estimator of $f$. 
 Note, however, that this approach requires analyticity of $\widehat{f}$ in a strip containing
 the imaginary axis, i.e., $f$ must have very light tails.
 %
In contrast, 
our construction does not require existence of $\widehat{f}(z)$ for \(z\) outside the imaginary axis; only the analyticity of 
$\widehat{g}(z)$ is required. 
\end{remark}

\section{Distribution of measurement errors}\label{sec:assumptions}

\subsection{Assumptions}\label{sec:gen-assumption}
Accuracy of the estimator $\tilde{f}_{s,h}(x_0)$ defined in (\ref{eq:f-sh}) 
will be studied under the following general assumptions
on the distribution of the measurement errors.
\begin{assumption}\label{as:general}
The Laplace transform $\widehat{g}(z)$ of the measurement error distribution 
exists in a vertical strip $\Sigma_g:=\{z\in \bC: \sigma_g^-<{\rm Re}(z)< \sigma_g^+\}$, $\sigma_g^-<0<\sigma_g^+$,  
and admits the following representation:
\begin{equation}\label{eq:factorization}
\widehat{g}(z)=  \frac{1}{\widehat{\psi}(z)} 
\prod_{k=1}^q \bigg(1- \frac{e^{a_k z}}{\lambda_k}\bigg)^{m_k},
\end{equation}
where  $a_1,\ldots,a_q$ are  positive real numbers, $|\lambda_k|=1$, $ k=1, \ldots, q,$  
$m_1,\ldots,m_q$ are non-negative integer numbers, and  
the pairs $(a_k, \lambda_k)$, $k=1,\ldots, q$  are distinct. 
The function  $\widehat{\psi}(z)$ is 
represented  as 
\begin{equation}\label{eq:psi}
  \widehat{\psi}(z)= 
  \widehat{\psi}_0(z)
  \prod_{k\in \Lambda} ( - a_k z)^{m_k} \prod_{k\not\in \Lambda} \bigg(1-\frac{1}{\lambda_k}\bigg)^{m_k}, 
 \end{equation}
 where $\Lambda:=\{k=1, \ldots, q: \lambda_k=1\}$,
 $\widehat{\psi}_0(z)$  is analytic, does not vanish in a vertical strip 
 $\Sigma_\psi:=\{z\in \bC: \sigma_\psi^-<{\rm Re}(z)<\sigma_\psi^+\}$  
 with  \mbox{\(\sigma_g^-\leq \sigma_\psi^-< 0 < \sigma_\psi^+\leq \sigma_g^+\)}, and $\widehat{\psi}_0(0)=1$.
\end{assumption}
Several remarks on Assumption~\ref{as:general} are in order. 
\begin{remark} 
{\rm (i)}. Assumption~\ref{as:general} states that  $\widehat{g}(z)$ factorizes into a product 
of two functions: 
the first function, $\prod_{k=1}^q(1-e^{a_k z}/\lambda_k)^{m_k}$, 
has zeros only on the imaginary
 axis, while the second one, $1/\widehat{\psi}(z)$, 
 does not have zeros 
 in $\Sigma_\psi\subseteq \Sigma_g$. The 
 latter fact follows from analyticity of $\widehat{\psi}(z)$ in 
 $\Sigma_\psi$.
 \par 
{\rm (ii)}.  Zeros of $\widehat{g}$ on the imaginary axis 
 are $z_{k,j}:=i({\rm arg}\{\lambda_k\}+2\pi j)/a_k$, $z_{k,j}\ne 0$, where $j\in \bZ$, $k=1, \ldots, q$, and 
 the multiplicity of each zero $z_{k,j}$  is equal to~$m_k$. 
 Assumption~\ref{as:general} implies that $\widehat{g}(z)$ has no zeros in $\Sigma_\psi\backslash \{z:
 {\rm Re}(z)=0\}$, and 
(\ref{eq:g-}) holds with  
 $S_g=\{z: {\rm Re}(z)\in (-\sigma_\psi^+, 0)\cup (0, -\sigma_\psi^-)\}$, that is, 
 $\kappa_g^-=-\sigma_\psi^+$ and $\kappa_g^+=-\sigma_\psi^-$.
 \par 
 {\rm (iii)}. The  form of $\widehat{\psi}(z)$ in (\ref{eq:psi}) 
 follows   from  (\ref{eq:factorization})
 and the fact that $\widehat{g}(0)=1$. 
 \end{remark} 
In addition to Assumption~\ref{as:general}
we require 
conditions on  the  growth of  
function $\widehat{\psi}(\cdot)$ in (\ref{eq:factorization})
on the imaginary axis.  These conditions 
are similar to the 
standard ones 
in  the {\em smooth case} [see condition~(B1) 
in Section~\ref{sec:introduction}]. 
 \begin{assumption}\label{as:psi}
Assume   
there exist constants $\omega_0>0$, $\gamma\geq 0$ and \mbox{$D_1>0$}, $D_2>0$  such that 
 \begin{align}\label{eq:smooth}
 &  D_1|\omega|^\gamma \leq  |\widehat{\psi}(i\omega)|\leq D_2|\omega|^\gamma,\;\;\;\forall |\omega|\geq \omega_0.
 \end{align}
In addition, suppose that 
\begin{align}\label{eq:smooth-1}
&  |\widehat{\psi}^{(j)}(i\omega)| \leq  D_{2+j} (1+|\omega|^{\gamma}),\;\;\;\forall \omega\in \bR 
\end{align}
for all natural $j\geq 1$ and some positive constants $D_3,D_4,\ldots.$ 
\end{assumption}
\begin{remark}
Condition (\ref{eq:smooth}), when imposed  on $1/\widehat{g}(i\omega)$, 
is rather standard in the literature; it corresponds to the so-called {\em smooth error densities}.
As for (\ref{eq:smooth-1}), similar  
restrictions on the derivatives of $\widehat{g}(i\omega)$
are usually imposed  in the proofs of lower bounds; see, e.g., 
\cite[Theorem~5]{Fan91}. We also note that all derivatives of the function $\widehat{\psi}(i\omega)$ exist 
due to analyticity of $\widehat{\psi}$. Furthermore, if the inequality
\[
\sup_{y\in [-r,r]} |\widehat{\psi} (i\omega+y)|\leq C(r)(1+|\omega|^\gamma), \quad \omega \in \mathbb{R} 
\]
holds for some  $r\in (0,\max\{-\sigma_\psi^-,\sigma_\psi^+\})$, and constant $C$ is independent of $\omega,$  
then the well-known Cauchy derivative estimates imply (\ref{eq:smooth-1}).
\end{remark}
 \subsection{Examples of distributions}\label{sec:examples}
 Assumptions~\ref{as:general} and~\ref{as:psi}
define a broad class of distributions containing
densities with characteristic functions that vanish on the real line. 
Moreover, discrete distributions are covered by Assumptions~\ref{as:general} and~\ref{as:psi}.
All this is illustrated in the  following examples.
\begin{example}[Uniform distribution] \label{ex:1}
{\rm 
Let $\epsilon~\sim U(-\theta, \theta),$ then 
\[
 \widehat{g}(z)= \frac{\sinh(\theta z)}{\theta z}=-\frac{e^{-\theta z}}{2\theta z}
 (1-e^{2\theta z}),\;\;z\in \bC.
\]
In this case 
representation (\ref{eq:factorization}) holds with 
$q=1$, \(m_1=1,\) $a_1=2\theta$, $\lambda_1=1$, and    
$\widehat{\psi}(z)=-2\theta z e^{\theta z}$,  $\widehat{\psi}_0(z)=e^{\theta z}$.
Clearly, 
$\widehat{\psi}$ satisfies Assumption~\ref{as:psi} with $\gamma=1$.
Note that $\widehat{g}(z)$ 
has simple zeros on the imaginary axis 
at $z_k=i\pi k/\theta$, $k=\pm 1, \pm 2, \ldots$, and 
$\Sigma_g=\Sigma_\psi=\bC$.
%
%
}
\end{example}
\begin{example}[Convolution of uniform distributions] \label{ex:2}
{\rm 
Consider a convolution of the uniform 
distributions $U(-\theta_k, \theta_k)$, $k=1,\ldots, q,$ with distinct parameters \(\theta_1,\ldots,\theta_q\), each of multiplicity \(m_k.\) In this case
\begin{eqnarray*}
 \widehat{g}(z)=\prod_{k=1}^q \left[\frac{\sinh(\theta_kz)}{\theta_k z}\right]^{m_k}
 =\frac{\exp\{-z \sum_{k=1}^q \theta_k m_k\}}{\prod_{k=1}^q (-2\theta_k z)^{m_k}} \prod_{k=1}^q 
 (1-e^{2\theta_kz})^{m_k}, \;\;\;z\in \bC.
\end{eqnarray*}
Therefore Assumption~\ref{as:general} holds with $a_k= 2\theta_k$,  $\lambda_k=1$,  
$k=1,\ldots, q$, 
$\Sigma_g=\Sigma_\psi=\bC$, and 
\begin{equation}\label{eq:psi-ex2}
\widehat{\psi}(z)=\prod_{k=1}^q (-2\theta_k z)^{m_k} e^{z\sum_{k=1}^q \theta_k m_k},\;\;\;\;
 \widehat{\psi}_0(z)= e^{z\sum_{k=1}^q \theta_k m_k}.
\end{equation}
Thus, $\widehat{\psi}(z)$ satisfies Assumption~\ref{as:psi} with $\gamma=m_1+\cdots+m_k$. 
Of special interest  is the case of
$m$ identical uniform distributions $U(-\theta, \theta)$. Here  
$q=1$, $\theta_1=\theta$, $m_1=m$, 
$\widehat{\psi}(z)=(-2\theta z)^m\exp\{m\theta z\}$,
$a_1=2\theta$ and $\lambda_1=1$. 
Note also that in this case $\gamma=m$.
}
\end{example}
\begin{example}[Discrete distributions] \label{ex:3} 
{\rm 
Let $\epsilon$ be a discrete random variable 
taking values in the set $\{0,\pm b,\ldots, \pm Mb\}$, $b>0,$  
with corresponding  probabilities $p_k$, $k=0,\pm 1, \ldots, \pm M$, where \mbox{$p_M\ne 0$}. Then 
\[
 \widehat{g}(z)= \sum_{k=-M}^M p_k e^{-b kz}= e^{-b Mz}\sum_{k=0}^{2M}
 p_{M-k} e^{b kz} = e^{-b M z} p_{M} P(e^{b z}), 
\]
where $P(x):=1+ \sum_{k=1}^{2M}(p_{M-k}/p_{M}) x^k$.
Let $\lambda_1, \ldots, \lambda_{2M}$ be the roots of the polynomial $P(x)$; 
 note that $\lambda_k\ne 1$, $\forall k=1, \ldots, 2M.$
Then we have 
\[
 \widehat{g}(z)=p_{M} e^{-b Mz} \prod_{k=1}^{2M} \Big(1- \frac{e^{b z}}{\lambda_k} \Big)
 = p_{M} e^{-b Mz} \prod_{k: |\lambda_k|\ne 1 } \Big(1- \frac{e^{b z}}{\lambda_k}\Big) 
 \prod_{k: |\lambda_k|=1 } \Big(1- \frac{e^{b z}}{\lambda_k}\Big), 
\]
and $\widehat{g}$ is an entire function, i.e., $\Sigma_g=\bC$.  
Representations (\ref{eq:factorization}) and \eqref{eq:psi} hold with  
\[
\widehat{\psi}(z)= \frac{e^{b M z}}{p_{M} \prod_{k: |\lambda_k|\ne 1}
(1- e^{b z}/\lambda_k)}, \quad \widehat{\psi}_0(z)= \frac{\widehat\psi(z)}{\prod_{k: |\lambda_k|= 1}
(1- 1/\lambda_k)}
\]
and $\Sigma_\psi=\big\{z: b^{-1}\ln (\lambda_-)< {\rm Re}(z) <b^{-1} \ln (\lambda_+)\big\}$, where
$\lambda_-:=\max\{ |\lambda_k|: |\lambda_k| < 1  \}$, and 
$\lambda_+:=\min\{ |\lambda_k|: |\lambda_k|>1\}$. 
In this example if all $\lambda_k$ with $|\lambda_k|=1$ are distinct, then 
 $q:=\#\{k: |\lambda_k|=1\}$, $a_1=\cdots=a_q=b$, and $m_1=\cdots=m_q=1$. 
It is obvious that Assumption~\ref{as:psi} holds with $\gamma=0$.
\par 
In the special case of  the Bernoulli distribution with the success probability  parameter 
$1/2$ we have  
$\widehat{g}(z)=\frac{1}{2}(1+e^z)$; hence (\ref{eq:factorization}) holds with 
$q=1$, $a_1=1$, $\lambda_1=-1$, $m_1=1$, and $\widehat{\psi}(z)=2$. If $\epsilon$ is a binomial random variable with the number of trials \(m\) and a success probability \(1/2\), then $\widehat{g}(z)=2^{-m}(1+e^z)^m$, and (\ref{eq:factorization})
holds with $q=1$, $a_1=1$, $\lambda_1=-1$, $m_1=m$, and $\widehat{\psi}(z)=2^m$.
}
\end{example}
\begin{example}[Convolution of uniform and smooth densities]\label{ex:4}
{\em 
Let $\varphi$ be a probability density with  Laplace transform $\widehat{\varphi}$ defined
in a strip $\Sigma_\varphi=\{z: \sigma_\varphi^-<{\rm Re}(z)<\sigma_\varphi^+\}$
such that $\widehat{\varphi}(z)\ne 0$, 
$\forall z\in \Sigma_\varphi$. Assume also that    
$|\widehat{\varphi}(i\omega)| \asymp |\omega|^{-\gamma}$ for some $\gamma>0$ as 
$|\omega|\to\infty$, that is, $\varphi$ 
satisfies the standard conditions of the smooth case.
%
Let $g$ be a convolution of the uniform density on $[-\theta, \theta]$ with $\varphi$; then 
\[
 \widehat{g}(z)=\frac{\sinh(\theta z)}{\theta z} \widehat{\varphi}(z)= -\frac{e^{-\theta z}\widehat{\varphi}(z)}{2\theta z} (1-e^{2\theta z}),  \;\;
 \sigma_{\varphi}^-<{\rm Re}(z) <\sigma_{\varphi}^+,
\]
and (\ref{eq:factorization})  obviously holds with $\widehat{\psi}(z)=-2\theta z e^{\theta z}/ \widehat{\varphi}(z)$. 
For instance, let $\varphi$ be 
a density of the Gamma distribution with parameters \(\gamma>0\) and \(\lambda>0\), that is, $\varphi(x)=\lambda^\gamma [\Gamma(\gamma)]^{-1} x^{\gamma-1} e^{-\lambda x}$, $x>0.$ Then  
$\widehat{\varphi}(z)= \lambda^\gamma (z+\lambda)^{-\gamma}$, ${\rm Re}(z)>-\lambda$, and $\widehat{\psi}(z)=-2\theta \lambda^{-\gamma} z e^{\theta z} (z+\lambda)^\gamma$.
} 
\end{example}
\section{Kernel representation}\label{sec:representation}
Under Assumption~\ref{as:general} kernel $L_{s,h}$  defined in (\ref{eq:L}) 
is rewritten as follows 
\begin{equation}\label{eq:L-general}
 L_{s,h}(t)=\frac{1}{2\pi}
 \int_{-\infty}^\infty  
 \frac{\widehat{K}((s+i\omega)h)\,\widehat{\psi}(-s-i\omega)}{\prod_{k=1}^q [1-e^{-a_k(s+i\omega)}/\lambda_k]^{m_k}}\,
 e^{(s+i\omega)t} \rd \omega,\quad s+i\omega\in S_g,
\end{equation}
where $S_g=\{z: {\rm Re}(z)\in (\kappa_g^-, 0)\cup (0, \kappa_g^+)\}$ is the set on which
$\widehat{g}(-z)$ does not vanish.  
Thus, for any $s\in (\kappa_g^-, 0)\cup (0, \kappa_g^+)$ the denominator of the integrand
in (\ref{eq:L-general}) is non--zero. Below 
we demonstrate  that $L_{s,h}$ can be  formally represented as an infinite series.
\subsection{Infinite series representation} 
To develop the infinite series representation 
we need the following notation.
According to Assumption~\ref{as:general}, 
the set of zeros of  $\widehat{g}(z)$ on the imaginary axis is 
determined by   three  
 $q$-tuples $a=(a_1,\ldots, a_q)$, $\lambda=(\lambda_1, \ldots, \lambda_q)$ and $m=(m_1,\ldots, m_q)$.
For a given vector $a=(a_1,\ldots, a_q)$ define  
\begin{equation*}
\sL:= \Big\{ a^Tj=\sum_{k=1}^q a_k j_k: j=(j_1,\ldots, j_q)\in \bZ_+^q\Big\}.     
\end{equation*} 
The set $\sL$ can be represented 
as an  ordered set of 
real numbers 
$\sL:=\{\ell_0, \ell_1, \ell_2,\ldots\}$,  where $0=\ell_0<\ell_1< \ell_2 < \ell_3< \cdots$. 
Define also 
\begin{equation}\label{eq:Rh}
 R_h(t) :=
  \frac{1}{2\pi} \int_{-\infty}^\infty 
  \widehat{K}(i\omega h) \widehat{\psi}(-i\omega) e^{i\omega t} \rd \omega,\;\;\;t\in \bR,
\end{equation}
and let
\[
 C_{j,m}:=\binom{j+m-1}{m-1}.
\]
In fact, $C_{j,m}$ is the number of {\em weak compositions of $j$ into $m$ parts} [see, e.g., 
\cite[p.~25]{Stanley}].
Recall  that an $m$--tuple $(i_1,\ldots, i_m)$ of non--negative integers 
with $i_1+\cdots+i_m=j$ is called a weak composition of~$j$ into $m$ parts. 
\begin{lemma}\label{lem:LL}
Let Assumption~\ref{as:general} hold, and 
 $\int_{-\infty}^\infty \big|\widehat{K}(i\omega h)\widehat{\psi}(-i\omega)\big| \rd \omega<\infty$.
\begin{itemize}
\item[{\rm (a)}] If $s\in (0, \kappa_g^+),$ then 
 \begin{align}
  & L_{s,h}(t) = L_{+,h}(t) :=  \sum_{\ell\in \sL} \cC^+_\ell R_h(t-\ell),
  \label{eq:C+}
  \\
  & \cC_\ell^+
  := 
\sum_{j: a^Tj=\ell}\;
\bigg[
\prod_{k=1}^q C_{j_k,m_k} 
\lambda_k^{-j_k}\bigg],
\;\;\ell\in \sL,
\label{eq:C+1}
\end{align}
provided that the summation  on
the right hand side of (\ref{eq:C+}) defines a  finite function for any~$t$. 
\item[{\rm (b)}] If $s\in (\kappa_g^-, 0)$ then 
\begin{align}
  & L_{s,h}(t) = L_{-,h}(t) :=  \sum_{\ell\in \sL} \cC^-_\ell R_h(t+a^Tm+\ell),
 \label{eq:C-}
 \\
  & \cC^-_\ell := 
\sum_{j: a^Tj=\ell}\;
\bigg[
\prod_{k=1}^q (-1)^{m_k} 
C_{j_k, m_k}
\lambda_k^{j_k+m_k}\bigg],
\;\;\ell\in \sL,
\label{eq:C-1}
\end{align}
provided that the summation on the right hand side of (\ref{eq:C-}) is finite for any $t$.
\end{itemize}
\end{lemma}
\begin{remark}\label{rem:2}
{\rm (i).} Lemma~\ref{lem:LL} shows that under Assumption~\ref{as:general}  
 kernel $L_{s,h}(t)$ can be represented as infinite linear combination of 
 one--sided translations of $R_h$, where the translation 
 parameter $\ell$ takes values in the set $\sL$. 
\par 
{\rm (ii).}
The coefficients $\{\cC_\ell^+\}$ and $\{\cC_\ell^-\}$ of the linear combination  are completely determined by the structure of the zero set of $\widehat{g}(z)$
on the imaginary axis. The sequences $\{\cC_\ell^+\}$, $\{\cC_\ell^-\}$ will play an 
important role in the sequel, and we call them {\em the zero set sequences}.
The definitions
in (\ref{eq:C+1}) and 
(\ref{eq:C-1}) imply that 
the coefficients $|\cC^+_\ell|$, $|\cC^-_\ell|$  may grow 
at most polynomially in $\ell$ as $\ell\to\infty$. Note also that $\cC_{\ell_0}^+=1$ and $\cC^-_{\ell_0}=\prod_{k=1}^q(-\lambda_k)^{m_k}$.
\end{remark}
\subsection{Kernel representation in specific problem instances} 
In general, determination of coefficients $\{\cC^+_\ell\}$ and $\{\cC^-_\ell\}$ 
in (\ref{eq:C+1}) and (\ref{eq:C-1}) is difficult. 
It is instructive to apply the result of Lemma~\ref{lem:LL} to some particular cases of 
Examples~\ref{ex:1}--\ref{ex:4}  
where 
the zero set sequences  
and the corresponding kernels can be explicitly calculated.
\paragraph{Uniform distribution}
This is the setting of Example~\ref{ex:1}.
Recall that here $q=1$,  $m_1=1$, $a_1=2\theta$, $\lambda_1=1$; hence $\sL=\{2\theta j: j=0, 1,\ldots\}$,
 and $\cC^+_\ell =1$, $\cC^-_\ell =-1$ for all $\ell=2\theta j$, $j=0,1,\ldots$. 
 Since $\widehat{\psi}(-i\omega)=2\theta i\omega e^{-i\theta \omega}$, 
 \[
  R_h(t)=\frac{1}{2\pi} \int_{-\infty}^\infty \widehat{K}(i\omega h) (2\theta i\omega) e^{i\omega(t-\theta)} \rd \omega = 
  \frac{2\theta}{h^2} K^\prime\Big(\frac{t-\theta}{h}\Big).
 \]
Thus, in view of (\ref{eq:C+}) and (\ref{eq:C-})
\begin{equation}\label{eq:L-uniform}
 L_{+, h}(t) = \frac{2\theta}{h^2}\sum_{j=0}^\infty K^\prime\Big(\frac{t-\theta(2j+1)}{h}\Big),\;
 L_{-,h}(t) = -\frac{2\theta}{h^2}\sum_{j=0}^\infty K^\prime
 \Big(\frac{t+\theta(2j+1)}{h}\Big).
\end{equation}
If $K$ is a bounded continuously differentiable kernel with finite support, and $h$ is small enough then formulas in (\ref{eq:L-uniform})
define   functions $L_{+,h}(t)$ and $L_{-, h}(t)$ which are finite for any~$t\in \bR$.  
%

\paragraph{Convolution of uniform distributions}
We consider two specific cases of Example~\ref{ex:2}.
\par\medskip 
(a). Consider convolution of $m$ identical uniform distributions $U(-\theta, \theta)$. In this case 
 $q=1$, $m_1=m$, 
$a_1=2\theta$ and $\lambda_1=1$. Thus, $\sL=\{2\theta j: j=0,1,\ldots\}$, $\cC^+_{2\theta j}= C_{j,m}$, 
$\cC^-_{2\theta j}=(-1)^mC_{j,m}$, $j=0,1,2,\ldots$. Since  $\widehat{\psi}(z)=(-2\theta z)^m\exp\{m\theta z\}$,
\[
 R_h(t)=\frac{1}{2\pi}\int_{-\infty}^\infty \widehat{K}(i\omega h) (2\theta i\omega)^m e^{i\omega (t-m\theta)}\rd \omega= 
 \frac{(2\theta)^m}{h^{m+1}}
 K^{(m)}\Big(\frac{t-m\theta}{h}\Big).
\]
Therefore
\begin{eqnarray*}
&& L_{+,h}(t)=\frac{(2\theta)^m}{h^{m+1}}\sum_{j=0}^\infty
 C_{j,m}K^{(m)}\Big(\frac{t-\theta(2j+m)}{h}\Big),\;
 \\
&& L_{-,h}(t)=\frac{(-2\theta)^m}{h^{m+1}}\sum_{j=0}^\infty
 C_{j,m}K^{(m)}\Big(\frac{t+\theta(2j+m)}{h}\Big).\;
\end{eqnarray*}
Similarly to the previous example,  if $K$ is $m$ times continuously differentiable with a finite support, and $h$ is small enough then 
the formulas define finite functions for any $t$.
\par\medskip 
(b).~Consider convolution of $q$ uniform distributions $U(-\theta_k, \theta_k)$, $k=1,\ldots, q$
with distinct $\theta_k$, $k=1,\ldots, q$.
In this case $a_k= 2\theta_k$, $\lambda_k=1$, $m_k=1$ for $k=1, \ldots, q$. Thus,
$\sL=\{ 2\sum_{k=1}^q \theta_k j_k: (j_1, \ldots, j_q)\in \bZ_+^q\}$, and if
$\ell= 2\sum_{k=1}^q \theta_k j_k^*$ for some $(j_1^*, \ldots, j_q^*)\in \bZ_+^q$ 
then $\cC^+_\ell$ is the number of non--negative
integer solutions $(x_1, \ldots, x_q)$ to the equation 
\[ 
\theta_1x_1+\cdots +\theta_q x_q= \theta_1j_1^*+\cdots+ \theta_q j^*_q,
\] 
and $\cC^-_\ell=(-1)^q\cC^+_\ell$.
It is clear that there is  at least one solution $(x_1,\ldots, x_q)=(j_1^*, \ldots, j_q^*)$; the total 
number of solutions
depends on $\theta_1, \ldots, \theta_q$. 
For instance,  assume that $\theta_k=r_k \theta_1$, $k=1,\ldots, q$,
where $1=r_1<r_2<\cdots<r_q$ are coprime integer numbers. 
Then $\cC^+_{\ell}$ with $\ell=2\theta_1\ell_*$ is the number of representations
of the integer number $\ell_*=j_1^*+ r_2j_2^*+\cdots +r_q j_q^*$ by non--negative integer linear combination of 
$r_1, \ldots, r_q$. 
By Schur's theorem [see, e.g., \cite[Section~3.15]{Wilf}] 
\begin{equation}\label{eq:Schur}
 \cC^+_\ell=\cC^+_{2\theta_1\ell_*} \sim \frac{\ell_*^{q-1}}{(q-1)! r_1r_2\cdots r_q},\;\;\; \ell_*\to\infty.
\end{equation}
It follows from (\ref{eq:psi-ex2}) that 
\[
\widehat{\psi}(-i\omega)= (2i\omega)^q \bigg(\prod_{k=1}^q \theta_k\bigg) \exp\{-i\omega \sum_{k=1}^q \theta_k\},
\]
and therefore 
\begin{eqnarray*}
 R_h(t)&=&\Big(\prod_{k=1}^q \theta_k\Big) 
 \frac{1}{2\pi} \int_{-\infty}^\infty \widehat{K}(i\omega h) (2i\omega)^q \exp\Big\{ i\omega\Big(t-\sum_{k=1}^q \theta_k\Big)\Big\}\rd \omega
 \\
 &=& \Big(2^q\prod_{k=1}^q \theta_k \Big) \frac{1}{h^{q+1}} K^{(q)}\Big(\frac{t-\sum_{k=1}^q \theta_k}{h}\Big). 
\end{eqnarray*}
Thus,
\begin{eqnarray*}
 L_{+,h}(t) &=&  \frac{\prod_{k=1}^q (2\theta_k)}{h^{q+1}} 
 \sum_{\ell\in \sL} \cC^+_\ell K^{(q)}\Big(\frac{t-\sum_{k=1}^q\theta_k -\ell}{h}\Big)
\\
L_{-,h}(t) &=&  \frac{\prod_{k=1}^q (2\theta_k)}{h^{q+1}} 
 \sum_{\ell\in \sL} \cC^-_\ell K^{(q)}\Big(\frac{t +\ell}{h}\Big),
 \end{eqnarray*}
 where $\cC^-_\ell=(-1)^q \cC^+_\ell$, and the sequence $\{\cC^+_\ell,\,\ell=2\theta_1j, \,j =0,1,\ldots\}$   
 satisfies
 (\ref{eq:Schur}).
 If kernel $K$ is $q$ times continuously differentiable and has bounded support then the last formulas define
 functions which are finite for any fixed~$t$.

\paragraph{Binomial distribution}
Assume that the measurement error distribution is  binomial with parameters $m$ and $1/2$; this is a 
particular case of Example~\ref{ex:3}. Here $q=1$, $a_1=1$, $\lambda_1=-1$, $m_1=m$ and $\widehat{\psi}(z)=2^m$. Hence 
$\sL=\bZ_+^1$, $\cC_\ell^+=C_{\ell, m}$, $\cC_\ell^-=(-1)^m C_{\ell, m}$, 
\[
 R_h(t) = \frac{1}{2\pi}\int_{-\infty}^\infty \widehat{K}(i\omega h)\widehat{\psi}(-i\omega) e^{i\omega t} \rd \omega = \frac{2^m}{h} K\Big(\frac{t}{h}\Big),
\]
and 
\begin{eqnarray*}
 L_{+, h}(t) = \frac{2^m}{h} \sum_{j=0}^\infty C_{j, m} K\Big(\frac{t-j}{h}\Big),\;
 L_{-, h}(t) = \frac{(-2)^m}{h} \sum_{j=0}^\infty C_{j, m} K\Big(\frac{t+j}{h}\Big)~.
\end{eqnarray*}


\section{Estimator and upper bounds on the risk}\label{sec:upper-bound}
Based on the general  ideas presented in Section~\ref{sec:idea}
and kernel representations developed in Section~\ref{sec:representation} 
we are now in a position to define the proposed estimator of $f$ and to study its accuracy.
\subsection{Estimator}
We assume that kernel $K$ is chosen to satisfy the following condition. 
\begin{itemize}
 \item[(K)] 
 Let $K\in C^\infty(\bR)$ be a function supported on $[-1,1]$ such that  
 for a fixed positive integer $k_0,$
 \[
  \int_{-1}^1 K(t) \rd t=1,\;\;\int_{-1}^1 t^j K(t) \rd t =0,\;\;j=1,2,\ldots, k_0.
 \]
\end{itemize}
Condition~(K) is standard in nonparametric kernel density estimation; clearly,   one can always construct kernel  
$K$  satisfying (K) with prescribed parameter $k_0$.
\par 
Let $N$ be a natural number, and  
denote 
\[
\sL_N:=\bigg\{a^Tj=\sum_{k=1}^qa_kj_k: j=(j_1,\ldots, j_q)\in \{0,1, \ldots, N\}^q\bigg\}. 
\]
The estimator of $f(x_0)$ is defined as follows
\begin{equation}\label{eq:estimator}
 \tilde{f}^{(N)}_{s,h}(x_0):= \frac{1}{n} \sum_{j=1}^n L_{s,h}^{(N)}(Y_j-x_0),\;\;s\in (\kappa_g^-,0)\cup 
 (0,\kappa_g^+),
\end{equation}
where we set 
\begin{equation}\label{eq:LL}
 L_{s,h}^{(N)}(t):= \left\{\begin{array}{ll}
                           L^{(N)}_{+,h}(t), & s\in (0, \kappa_g^+),\\*[2mm]
                           L^{(N)}_{-,h}(t), & s\in (-\kappa_g^-, 0),
                          \end{array} \right. 
\end{equation}
and 
\begin{equation}\label{eq:L+-}
 L_{+,h}^{(N)}(t):= \sum_{\ell\in \sL_N} \cC^+_\ell R_h(t-\ell),\;\;
 L_{-,h}^{(N)}(t):= \sum_{\ell\in \sL_N} \cC^-_\ell R_h(t+\ell).
\end{equation}
In what follows we will write $\tilde{f}_{+, h}^{(N)}(x_0)$ and $\tilde{f}_{-,h}^{(N)}(x_0)$ for the estimator
(\ref{eq:estimator}) associated with $s\in (0, \kappa_g^+)$ and $s\in (\kappa_g^-, 0),$ respectively.
Let finally
\begin{equation}\label{eq:estimator-1}
 \tilde{f}_h^{(N)}(x_0):= \left\{\begin{array}{ll}
                                  \tilde{f}_{+, h}^{(N)}(x_0), & x_0\geq 0,\\*[2mm]
                                  \tilde{f}_{-,h}^{(N)}(x_0), & x_0<0.
                                 \end{array}
\right.
\end{equation}
\par 
Recall that the function $R_h$ and the sequences $\{\cC^+_\ell, \ell\in \sL\}$, $\{\cC^-_\ell, \ell\in \sL\}$
are defined in (\ref{eq:Rh}),  (\ref{eq:C+}) and (\ref{eq:C-}), respectively.
The estimator construction follows {\em the linear functional strategy} of Section~\ref{sec:idea}
in conjunction with the kernel representation developed in Section~\ref{sec:representation}. 
Note that we truncate the infinite series kernel representation by the cut--off parameter~$N$; 
this introduces some bias
but ensures that the integral on the left hand side of (\ref{eq:L-K}) is absolutely convergent.  
The estimator $\tilde{f}^{(N)}_{s,h}(x_0)$ requires specification of   $h$ and $N$; this 
will be done in the sequel. 
\par

\subsection{Functional classes}
Now we 
define functional classes over which accuracy of the proposed estimators
will be assessed. The next two definitions introduce standard classes of smooth functions.
\par 
\begin{definition}
\label{def:1}
 Let $A>0$, $\alpha>0$ be real numbers. 
  We say that  a probability density $f$ belongs to the functional class $\sH_\alpha(A)$ if 
  $f$ is $\lfloor \alpha\rfloor := \max\{k\in \bN\cup \{0\}:  k<\alpha\}$ 
  times continuously differentiable and 
  \[
   |f^{(\lfloor \alpha \rfloor)}(x)- f^{(\lfloor \alpha \rfloor)}(x^\prime)| \leq A |x-x^\prime|^{\alpha-\lfloor \alpha\rfloor},
   \;\;\;\forall x, x^\prime \in \bR.
  \]
 \end{definition}
\begin{definition}
For real numbers $A>0$ and $\alpha>1/2$ 
  we say that a probability density $f$ belongs to the functional class $\sS_\alpha(A)$ if 
   $f$ is $\lfloor \alpha\rfloor := \max\{k\in \bN\cup \{0\}:  k<\alpha\}$ 
  times differentiable and 
  \[
   \big\|f^{(\lfloor \alpha \rfloor)}(\cdot+ t)- f^{(\lfloor \alpha \rfloor)}(\cdot)\big\|_2 \leq A |t|^{\alpha-\lfloor \alpha\rfloor},
   \;\;\;\forall t \in \bR.
  \]
\end{definition}
\par 
We will also consider classes of probability densities with bounded moments. 
\begin{definition}
Let $p>0$ and $B>0$ be real numbers. 
We say that a probability density $f$ belongs to the functional class 
$\sM_p(B)$ if 
 \[
 \max\Big\{\|f\|_\infty,\; \max_{0<r\leq p}\int_{-\infty}^\infty |t|^r f(t)\rd t\Big\}\leq B.
 \]
 We also denote $\sM_p^\prime (B)$ the class of all densities $f$ from $\sM_p(B)$ satisfying the following additional condition: 
\begin{equation}
\max_{0\leq j\leq  p}\int_{-\infty}^\infty \frac{|\widehat{f}^{(j)}(i\omega)|}{1+ |\omega|^\gamma} \;\rd \omega \leq B,
\label{eq:the-condition-2}
\end{equation}
where $\gamma\geq 0$ is a constant appearing in Assumption~\ref{as:psi}. 
\end{definition}

\begin{remark}\label{rem:class}
{\rm (i).} 
 Condition (\ref{eq:the-condition-2}) is rather mild. Note that  $f\in \sM_p(B)$ implies
boundedness of  $|\widehat{f}^{(j)}(i\omega)|$ for all $j=0,\ldots, p$, and in (\ref{eq:the-condition-2}) we require
integrability of $\widehat{f}^{(j)}(i\omega)$ with the weight $(1+|\omega|^\gamma)^{-1}$.
If $\gamma>1$ and $f\in \sM_p(B_1),$ then (\ref{eq:the-condition-2}) holds trivially:
$f\in \sM_p^\prime (B_2)$ with $B_2=c_\gamma B_1$ where 
$c_\gamma := \int_{-\infty}^\infty (1+|\omega|^\gamma)^{-1}\rd \omega$. Therefore in the definition of 
$\sM^\prime_p(B)$ the restriction (\ref{eq:the-condition-2}) is active 
only if $\gamma\leq 1$.
\par
{\rm (ii).}  
If $f\in \sH_\alpha(A),$ then $f$ is uniformly bounded above by a constant  depending on $A$ only. 
However, for the sake of convenience, we explicitly require boundedness of $f$ in the definition of the class $\sM_p(B)$.  
 \end{remark} 
We also denote 
\begin{align*} 
 &\sF_{\alpha,  p}(A, B)  := \sH_{\alpha}(A) \cap \sM_{p} (B),\;\;
 \sF^\prime_{\alpha, p}(A, B)  := \sH_\alpha(A)\cap \sM^\prime_p(B),
 \\
 &\sG_{\alpha, p}(A,B):=\sS_{\alpha}(A) \cap \sM_{p} (B),\;\;
 \sG^\prime_{\alpha, p}(A, B)  := \sS_\alpha(A)\cap \sM^\prime_p(B).
\end{align*} 
 \subsection{Upper bounds}
In this section we  derive 
an  upper bound on the maximal risk of the estimator (\ref{eq:estimator}) 
under Assumptions~\ref{as:general}, \ref{as:psi}, and under the following additional condition on the growth of zero set sequences $\{\cC_\ell^+\}$ and~$\{\cC_\ell^-\}$.
\begin{assumption}\label{as:C+-}
 Assume that 
 \begin{equation*}
  \sum_{\ell\in \sL\backslash \{0\}} \max\{|\cC_\ell^+|, |\cC_\ell^-|\} \,\ell^{-\nu} \leq C_0 <\infty
 \end{equation*}
for some $C_0>0$ and $\nu> 1.$  
\end{assumption}
\par

\begin{theorem}\label{th:upper}
Suppose that Assumptions~\ref{as:general}, \ref{as:psi} and~\ref{as:C+-} hold.
Let $\tilde{f}_h^{(N)}(x_0)$ be associated with kernel $K$ satisfying 
 the condition (K) with  $k_0\geq \alpha + 1$.
\begin{itemize}
\item [{\rm (a)}] 
Assume that $f\in \sF^\prime_{\alpha, p}(A, B)$ with 
$p\geq 2\nu.$
Let 
$h=h_*:=\big[B(A^2 n)^{-1}\big]^{1/(2\alpha+2\gamma+1)}$ and 
$N\geq \big(A^{-2\gamma+1} B^{2\gamma+\alpha} 
  n^{\alpha+1}\big)^{1/p(2\alpha+2\gamma+1)}$.
Then for all large enough $n$ one has
\[
 \cR_{n, \Delta_{x_0}}\big[\tilde{f}_{h_*}^{(N)}; \sF^\prime_{\alpha, p}(A,B)\big] \leq C_1 A^{\frac{2\gamma+1}{2\alpha+2\gamma+1}} \big(Bn^{-1})^{\frac{\alpha}{2\alpha+2\gamma+1}},
\]
where $C_1$ may depend on $\alpha,$  and $p$  only.
\item[{\rm (b)}] 
Let $f\in \sG^\prime_{\alpha, p}(A, B)$ 
with  
$p\geq 2\nu+1.$
Let 
$N\geq \big(A^{-2\gamma+1} B^{2\gamma+\alpha} 
  n^{\alpha+1}\big)^{2/(2p-1)(2\alpha+2\gamma+1)}$
and $h=h_*$. Then for all large enough $n$ one has 
\[
 \cR_{n, \Delta_{2}}\big[\tilde{f}_{h_*}^{(N)}; \sG^\prime_{\alpha, p}(A,B)\big] \leq C_2 A^{\frac{2\gamma+1}{2\alpha+2\gamma+1}} \big(Bn^{-1})^{\frac{\alpha}{2\alpha+2\gamma+1}},
\]
where $C_2$ may depend on $\alpha$  and $p$  only.
\end{itemize}
\end{theorem}
\begin{remark}
{\rm (i).} It is well known that under 
assumptions $\widehat{g}(i\omega) \ne 0$, \(\omega \in \bR\) and $\widehat{g}(i\omega)\asymp |\omega|^{-\gamma}$,
\mbox{$|\omega|\to\infty$}, 
we have  $\cR_{n, \Delta_{x_0}}^*[\sH_\alpha(A)]\asymp n^{-\alpha/(2\alpha+2\gamma+1)}$  and
 $\cR_{n, \Delta_{2}}^*[\sS_\alpha(A)]\asymp n^{-\alpha/(2\alpha+2\gamma+1)}$ as \mbox{$n\to\infty$}.
Theorem~\ref{th:upper} provides conditions on $f$ that guarantee 
 the standard  rate of convergence 
 in the case  when \(\widehat{g}(i\omega)\) may have zeros. In particular, 
conditions 
$f\in \sM^\prime_{p}(B)$  with  $p\geq 2\nu$ and $p\geq 2\nu+1$ 
are sufficient  in order to ensure
 the standard rate $n^{-\alpha/(2\alpha+2\gamma+1)}$ for the pointwise and $\bL_2$--risks, respectively.
It is worth noting that the $\bL_2$--risk bound  requires stronger condition; 
as we will show in Section~\ref{sec:lower-bound}, 
this is  an intrinsic feature 
of the problem. 
 \par {\rm (ii).}
 The result of Theorem~\ref{th:upper} is rather general: it holds for any configuration 
 of zeros of $\widehat{g}$ and  function $\widehat{\psi}$ satisfying Assumption~\ref{as:psi}. One interesting implication
 is that for discrete error distributions such as in Example~\ref{ex:3}, the achievable rate of convergence 
 is $n^{-\alpha/(2\alpha+1)},$ provided that $f$ has a finite absolute moment of  sufficiently large order. Note 
 that $n^{-\alpha/(2\alpha+1)}$
 is the minimax rate of convergence 
 in the problem of density estimation  from direct  
 i.i.d. observations.%
 \par 
{\rm (iii).} If the measurement error distribution is uniform, then Assumption~\ref{as:C+-}
holds with $\nu=1+\varepsilon$ for arbitrary small $\varepsilon >0$. Therefore 
Theorem~\ref{th:upper} implies that
the standard rate of convergence $n^{-\alpha/(2\alpha+3)}$ 
of the pointwise and $\bL_2$--risks is achieved  if $f\in \sM_p^\prime(B)$ with $p> 2$ and  $p>3$, respectively.
 \end{remark}
\par 
In some specific cases when closed form expressions for the zero sequences 
$\{\cC^+_\ell\}$ and $\{\cC_\ell^-\}$ and function $\widehat{\psi}(z)$  are available,  
the conditions of 
Theorem~\ref{th:upper}
can be  relaxed.
We demonstrate this in the next section.

\section{Specific problem instances}\label{sec:special}

In this section we consider specific distributions of measurement errors for which  
conditions of Theorem~\ref{th:upper} can be relaxed. 
\subsection{Convolution of uniform distributions}
Consider a particular case of Example~\ref{ex:2} where   $\widehat{g}(z)=[\sinh(\theta z)/(\theta z)]^m$, $m\geq 1$.
This setting also covers Example~\ref{ex:1} that corresponds to  $m=1$. 
Recall that in this case $\cC^+_{2\theta j}=C_{j,m}$, $\cC^-_{2\theta j}=(-1)^m C_{j,m}$, and therefore
Assumption~\ref{as:C+-} is valid for any $\nu>m$. Then Theorem~\ref{th:upper} implies that
the pointwise and $\bL_2$--risks converge to zero at the standard rate 
provided that  $f\in \sM^\prime_p(B)$ with $p>2m$, and 
$p>2m+1$, respectively.
In fact, as the following result demonstrates,
these conditions are too strong: 
the standard
rate is achievable if 
$p\geq 2m-2$ for the pointwise risk, and if $p> 2m-1$ for the $\bL_2$--risk. 
Recall that 
$\widehat{\psi}(z)=(-2\theta z)^m \exp\{m\theta z\}$ and  
\[
R_h(t)=
 \frac{(2\theta)^m}{h^{m+1}}
 K^{(m)}\Big(\frac{t-m\theta}{h}\Big).
 \]
 The corresponding kernels are   
 \begin{eqnarray}\label{eq:L+conv}
&& L^{(N)}_{+,h}(t)=\frac{(2\theta)^m}{h^{m+1}}\sum_{j=0}^N
 C_{j,m}K^{(m)}\Big(\frac{t-\theta(2j+m)}{h}\Big),\;
 \\
&& L^{(N)}_{-,h}(t)=\frac{(-2\theta)^m}{h^{m+1}}\sum_{j=0}^N
 C_{j,m}K^{(m)}\Big(\frac{t+\theta(2j+m)}{h}\Big).\;
\label{eq:L-conv}
 \end{eqnarray}
\begin{theorem}\label{th:conv-uniform}
Let  $\widehat{g}(z)=[\sinh(\theta z)/(\theta z)]^m$, $m\in \mathbb{N}$. 
Let $K$ be a kernel satisfying 
 the condition (K) with  $k_0\geq \alpha + 1$, and 
let $\tilde{f}_h^{(N)}$ denote the estimator defined in 
(\ref{eq:estimator})--(\ref{eq:estimator-1}) 
 and associated with the kernels  
(\ref{eq:L+conv}) and (\ref{eq:L-conv}).
\begin{itemize}
 \item[{\rm (a)}] 
 Assume that $f\in \sF_{\alpha, p}(A, B)$ with $p\geq 2m-2$ if $m>1$, 
 and   $p>0$ if $m=1$.  
 Let
 $h=h_*:=
 \big[B
 (A^2 n)^{-1}\big]^{1/(2\alpha+2m+1)}$, and 
$N\geq  \big(A^{-2m+1} B^{2m+\alpha} 
  n^{\alpha+1}\big)^{1/p(2\alpha+2m+1)}$.
Then for large enough $n,$
\[
 \cR_{n, \Delta_{x_0}}[\tilde{f}_{h_*}^{(N)}; \sF_{\alpha, p}(A, B)] \leq C_1 A^{\frac{2m+1}{2\alpha+2m+1}} \big(B
 n^{-1})^{\frac{\alpha}{2\alpha+2m+1}},
\]
where $C_1$ may depend on $\alpha$, $m$ and $\theta$ only.
\item[{\rm (b)}] 
Let $f\in \sG_{\alpha, p}(A, B)$ with  $p> 2m-1$, 
$N\geq (B^{2m+\alpha} A^{-2m+1} n^{\alpha+1})^{2/(2p-1)(2\alpha+2m+1)}$, and $h=h_*$.
Then for all large enough $n$
\[
 \cR_{n, \Delta_2}[\tilde{f}_{h_*}^{(N)}; \sG_{\alpha, p}(A, B)] \leq C_2 A^{\frac{2m+1}{2\alpha+2m+1}} \big(B n^{-1})^{\frac{\alpha}{2\alpha+2m+1}},
\]
where $C_2$ may depend on $\alpha$, $m$ and $\theta$  only.
\end{itemize}
\end{theorem}
\begin{remark}
{\rm (i).} In contrast to the proof of Theorem~\ref{th:upper},  
the proof of Theorem~\ref{th:conv-uniform} relies on a closed form expressions for 
kernels $L_{\pm,h}^{(N)}(t)$ [cf. (\ref{eq:L+conv}), (\ref{eq:L-conv})].
In this case   support of function $R_h$ has a ``small'' length, and  this fact is 
crucial for  relaxing assumptions of Theorem~\ref{th:upper}.
 \par 
 {\rm (ii).} Theorem~\ref{th:conv-uniform} shows that the standard rate of 
 convergence  \(n^{-\alpha/(2\alpha+2\gamma+1)}\) is achieved   
 \begin{itemize}
 \item 
 by the maximal pointwise risk over $\sH_\alpha(A),$ if
 $f\in \sM_{2m-2}(B)$ when $m>1$, and $f\in \sM_{\varepsilon}(B)$ for any
 $\varepsilon>0$ when $m=1$; 
 \item by the maximal $\bL_2$--risk over $\sS_\alpha(A),$ if 
 $f\in \sM_{2m-1+\epsilon}(B)$, $\epsilon>0$.
  \end{itemize}
 In contrast, Theorem~\ref{th:upper} requires  $f\in \sM_{p}(B)$ with $p>2m$  and 
 $f\in \sM_{p}(B)$, $p>2m+1,$ respectively. The difference in 
 conditions  is particularly noticeable  in the case of the uniform 
 density where $m=1$. Indeed, while Theorem~\ref{th:upper} requires finiteness of $p$-th moment with $p>2$ for the pointwise risk and $p>3$ for the 
 $\bL_2$--risk, Theorem~\ref{th:conv-uniform}
 shows that it is sufficient to require 
 conditions $f\in \sM_\varepsilon(B)$ and $f\in \sM_{1+\varepsilon}(B)$,
 $\varepsilon>0$.
It also follows from the proof that 
 assumption $f\in \sM_\varepsilon(B)$ for the pointwise risk
can be further relaxed: any uniform decrease of $f(x)$ as $|x|\to\infty$ will be sufficient. 
 \end{remark}

\subsection{Binomial distribution}
In this section we consider a specific case of 
Example~\ref{ex:3} of Section~\ref{sec:examples}, 
where the measurement errors have binomial distribution 
with parameters $m$ and $1/2$. Here 
$\widehat{g}(z)= 2^{-m} (1+ e^z)^m$.
Recall that in this case $\sL=\bZ_+^1$, $\cC_\ell^+=C_{\ell, m}$, $\cC^-_\ell=(-1)^m C_{\ell, m}$, 
\[
 R_h(t) = \frac{1}{2\pi}\int_{-\infty}^\infty \widehat{K}(i\omega h)\widehat{\psi}(-i\omega) e^{i\omega t} \rd \omega = \frac{2^m}{h} K\Big(\frac{t}{h}\Big)
\]
and 
\begin{equation}\label{eq:L-binom}
 L^{(N)}_{+, h}(t) = \frac{2^m}{h} \sum_{j=0}^N C_{j, m} K\Big(\frac{t-j}{h}\Big),\;\;
 L^{(N)}_{-, h}(t) = \frac{(-2)^m}{h} \sum_{j=0}^N C_{j, m} K\Big(\frac{t+j}{h}\Big).
\end{equation}

\begin{theorem}\label{th:binomial}
Let  $\widehat{g}(z)= 2^{-m} (1+ e^z)^m$, $m\geq 1$. Fix some \(\alpha>0.\) Let $K$ be a kernel  satisfying 
 condition (K) with  $k_0\geq \alpha + 1$, and 
let $\tilde{f}_h^{(N)}$ denote the estimator defined in 
 and (\ref{eq:estimator})--(\ref{eq:estimator-1}) 
 and associated with the kernels  in (\ref{eq:L-binom}).

\begin{itemize}
 \item[{\rm (a)}] Assume that $f\in \sF_{\alpha, p}(A, B)$ with 
 $p\geq 2m-2$ if $m>1$, and  $p>0$ if $m=1$. 
 Let 
 $h=h_*:=[B(A^2 n)^{-1}]^{1/(2\alpha+1)}$, and 
$N\geq (A B^{\alpha}n^{\alpha+1})^{1/p(2\alpha+1)}$.
 Then for large enough~$n,$
\[
 \cR_{n, \Delta_{x_0}} [\tilde{f}_{h_*}^{(N)}; \sF_{\alpha, p}(A, B)] \leq C_1 A^{\frac{1}{2\alpha+1}} \big(Bn^{-1})^{\frac{\alpha}{2\alpha+1}},
\]
where $C_1$ may depend on $\alpha$  only.
\item[{\rm (b)}] 
Assume that $f\in \sG_{\alpha, p}(A, B)$ for some 
$p> 2m-1$. 
Let $N\geq (A B^{\alpha} n^{\alpha+1})^{2/(2p-1)(2\alpha+1)}$, and $h=h_*$.
 Then for all large enough $n,$
\[
 \cR_{n, \Delta_{2}} [\tilde{f}_{h_*}^{(N)}; \sG_{\alpha, p}(A, B)] \leq C_2 A^{\frac{1}{2\alpha+1}} \big(Bn^{-1})^{\frac{\alpha}{2\alpha+1}},
\]
where $C_2$ may depend on $\alpha$  only.
\end{itemize}
\end{theorem}
The proof is omitted as it goes along the same lines as the proof of Theorem~\ref{th:conv-uniform} with minor modifications.
\section{Lower bounds: necessity of moment conditions}\label{sec:lower-bound}
Theorem~\ref{th:conv-uniform} shows that 
if the error distribution is the $m$--fold convolution of the uniform distributions on $[-\theta, \theta],$ then the maximal pointwise and $\bL_2$--risks on 
the classes $\sH_\alpha(A)$ and $\sS_\alpha(A)$ 
converge to zero at 
the standard  rate $n^{-\alpha/(2\alpha+2m+1)}$, 
provided that moment conditions hold.
 The following theorem demonstrates that
these moment conditions are also necessary.
\begin{theorem}\label{th:lower-bound-2}
 Let $\widehat{g}(i\omega)=[\sin(\theta\omega)/(\theta\omega)]^m$, where $m\geq 1$ 
 is integer.
\begin{itemize}
 \item[{\rm (a)}]
 If  $f\in \sF_{\alpha, p}(A, B)$ with  
 $0 < p < 2m-2$
 then 
\begin{equation}\label{eq:low-b-1}
 \liminf_{n\to\infty} \Big\{ n^{\alpha/(2\alpha+2m+1)}\, \cR_{n,\Delta_{x_0}}^*[\sF_{\alpha, p}(A,B)]\Big\} =\infty.
\end{equation}
\item[{\rm (b)}]
Moreover, if $f\in \sG_{\alpha, p}(A, B)$ with 
 $p<2m-1$
then 
\begin{equation}\label{eq:low-b-2}
 \liminf_{n\to\infty} \Big\{ n^{\alpha/(2\alpha+2m+1)}\, \cR_{n, \Delta_2}^*[\sG_{\alpha, p}(A,B)]\Big\} =\infty.
\end{equation}
\end{itemize}
\end{theorem}
 \begin{remark}
  The theorem states that under conditions 
  $0<p<2m-2$ and $p<2m-1$ 
  the standard rates of convergence cannot be attained in estimation with pointwise and 
  $\bL_2$--risks,
  respectively.
  On the other hand, we have constructed an estimator that
  achieves the standard rate of convergence, provided that 
  $p\geq 2m-2$, $m>1$ and $p>0$, $m=1$ in the former case,
  and $p> 2m-1$ in the latter case; 
  see  Theorem~\ref{th:conv-uniform}. Thus, the indicated  
  moment conditions are necessary 
  for convergence of the risks at the standard rate.
\end{remark}

\section{Concluding remarks}\label{sec:conclusion}
We close this paper with some concluding remarks.
\par\smallskip
1.~The proposed estimator $\tilde{f}_h^{(N)}(x_0)$ in (\ref{eq:estimator-1}) is 
associated with the one--sided kernels $L_{+, h}$ and $L_{-,h}$, 
which were used for positive and negative  values of $x_0,$ respectively. 
This  definition was adopted for the sake of convenience and unification of proofs. 
In fact, a closer inspection of the proofs 
of Theorems~\ref{th:upper} and \ref{th:conv-uniform} shows that
for any $x_0$ one can construct an estimator relying on any one of these two kernels
with the same risk guarantees. In this case 
the parameter $N$ should be chosen depending on $x_0$.
\par \smallskip 
2. In this paper we considered the functional class $\sM_p(B)$ of densities satisfying 
certain moment conditions. 
It is worth noting that the proposed estimators can be analyzed under other assumptions as well. 
For instance, if the support of $f$
has a finite left endpoint, then  there is no need to assume that $f\in \sM_p(B)$. 
Indeed, the proof
of Theorem~\ref{th:conv-uniform} shows that the accuracy of 
$\tilde{f}_{+, h}^{(N)}(x_0)$  ($\tilde{f}_{+, h}^{(N)}(x_0)$) is  determined by the right (left) tail of $f$.
Therefore if the support of $f$ has a  finite left endpoint, then it is reasonable to use 
the estimator $\tilde{f}_{-, h}^{(N)}(x_0)$
whose risk will converge to zero at the standard rate. 
This fact connects our result to those of 
\cite{GroenJong03} and \cite{Delaigle-Meister}.
\par\smallskip 
3. The following lower bounds on the minimax risks
$\cR_{n, \Delta_{x_0}}^*[\sH_\alpha(A)]$ and $\cR_{n, \Delta_2}^*[\sS_\alpha(A)]$ 
can be extracted from the proof of Theorem~\ref{th:lower-bound-2}. 
If the measurement error distribution is 
a  $m$--fold  convolution of the uniform distribution
then 
for any~$\delta>0$
\begin{align*}
 & \cR_{n, \Delta_{x_0}}^*[\sH_{\alpha}(A)] \geq C_1(\delta) \phi_{n, \Delta_{x_0}}(m),\;\;
 \phi_{n, \Delta_{x_0}}(m):= \bigg(\frac{A^{2m/\alpha}}{n}\bigg)^{\frac{\alpha}{2m\alpha+2m+1}-\delta},
 \\
 & \cR_{n, \Delta_2}^*[\sS_\alpha(A)] \geq C_2(\delta) \phi_{n, \Delta_2}(m),\;\;\phi_{n, \Delta_2}(m):=
 \bigg(\frac{A^{2m/\alpha}}{n}\bigg)^{\frac{\alpha}{4m\alpha+2}-\delta}.
\end{align*}
Observe  that 
$\phi_{n, \Delta_{x_0}}(m)\gg n^{-\alpha/(2\alpha+2m+1)}$ if $m>1$, and 
$\phi_{n, \Delta_{2}}(m)\gg n^{-\alpha/(2\alpha+2m+1)}$  for all $m \geq 1$.
These results should be compared with the upper bounds in
Theorem~\ref{th:conv-uniform}. In particular, even in the case of the uniform error density  
there is a significant difference in 
the behavior of the minimax risks 
$\cR^*_{n, \Delta_2}[\sS_\alpha(A)\cap \sM_{1+\varepsilon}(B)]$ and $\cR^*_{n,\Delta_2}[\sS_\alpha(A)]$:  
while the former is of the order  $n^{-\alpha/(2\alpha+3)}$, the latter one converges to zero at the rate slower than 
$n^{-\alpha/(4\alpha+2)-\delta}$
for any small $\delta>0$.
It is worth noting that 
some lower bounds on the minimax $\bL_2$--risk are reported in \cite{Hall-Meister} and \cite{Meister-1}. However, these bounds are not directly comparable with ours
since the considered functional classes and  assumptions on $\widehat{g}(z)$ in the above papers are different from ones adopted in our paper.  Moreover we mainly focus here on the minimal conditions needed to preserve the  standard convergence rates and do not consider the problem of constructing optimal  (in minimax sense) estimators in the case where these conditions are violated. 
\par\medskip 
4. We focused on the setting when characteristic function of measurement errors
has zeros on the imaginary axis and decreases at a polynomial rate. 
This corresponds to the case of {\em smooth error densities}. 
The {\em super--smooth case}   when the characteristic function 
decreases at an exponential rate can be also considered within the proposed framework.
This assumption leads to slow logarithmic rates, and it can be shown that
zeros of the error characteristic function do not affect the minimax rates of convergence, i.e.,
the standard minimax rates are preserved with no additional tail conditions.
\par\smallskip 
5. The proposed estimators are not adaptive in the sense that they require
knowledge  of the underlying functional classes. However, based on the results of this paper, adaptive estimators
can be developed using standard methods [see, e.g., \cite{Lepski90} and 
\cite{GL11}].  We do not pursue this direction in the current paper.
\par\smallskip 
6.  \cite{Johnstone-Raimondo} considered a closely related 
problem of signal deconvolution in the periodic  Gaussian
white noise model  
$\rd y(t)= (f\ast g)(t)\rd t +\sigma \rd W(t)$, $t\in [-1, 1]$, where 
$g$ is a boxcar kernel, that is, $g(t)=(2\theta)^{-1}{\bf 1}_{[-\theta, \theta]}(t)$, and $\{W(t)\}$ is the standard
two--sided Wiener process. If $\theta$ is a rational number then the signal $f$ 
is non--identifiable. Assuming that $\theta$ is irrational, 
\cite{Johnstone-Raimondo} studied 
 behavior of the minimax $\bL_2$--risk over the classes of 
 ellipsoids and hyperrectangles defined on  
the Fourier coefficients of $f$. They show that 
the minimax rates of convergence for the $\bL_2$--risk are affected by an oscillating 
behavior of the Fourier coefficients of the boxcar kernel.
Our results suggest that   if the assumption of periodicity of $f$ and $g$ is dropped then 
the minimax $\bL_2$--risk over the class  $\sS_\alpha(A)$  should be  of
the standard order
$(\sigma^2)^{-\alpha/(2\alpha+3)}$. We plan to study these signal deconvolution models 
in our future research.

\appendix 
\section*{Appendix}\label{sec:proofs}

\subsection*{Proof of Lemma~\ref{lem:K-L}}
Fix $s\in (\kappa_g^-, 0)\cup (0, \kappa_g^+)$. 
By the Fubini's theorem
\begin{eqnarray*}
 \int_{-\infty}^\infty L_{s,h}(y-x_0) f_Y(y)\rd y = 
 \int_{-\infty}^\infty L_{s,h}(y-x_0) \int_{-\infty}^\infty  f(x) \rd G(y-x) \rd x 
\\
=\int_{-\infty}^\infty f(x) \int_{-\infty}^\infty L_{s,h}(y-x_0) \rd G(y-x)  \rd x.
 \end{eqnarray*}
Now we show that for almost all $x,$ 
\begin{equation}\label{eq:L-g-K}
 \int_{-\infty}^\infty L_{s,h}(y-x_0) \rd G(y-x)  = \frac{1}{h} K\bigg(\frac{x-x_0}{h}\bigg).
\end{equation}
Applying the bilateral Laplace transform to the left hand side of the previous display formula we 
obtain 
\begin{eqnarray*}
 \int_{-\infty}^\infty e^{-zx} \int_{-\infty}^\infty L_{s,h}(y-x_0) \rd G(y-x)  \rd x = 
 \widehat{g}(-z) \widehat{L}_{h}(z) e^{-zx_0}.
\end{eqnarray*}
In view of (\ref{eq:L-hat}), the function on the right hand side 
is analytic  and  equal to
$\widehat{K}(zh) e^{-zx_0}$ on $S_g$. On the other hand, 
\[
 \int_{-\infty}^\infty e^{-zx} \frac{1}{h} K\bigg(\frac{x-x_0}{h}\bigg) \rd x = e^{-z x_0} \widehat{K}(zh),\;\;\;
 z\in \bC.
\]
Thus, the bilateral Laplace transforms of the functions on both sides
of (\ref{eq:L-g-K}) coincide on~$\bC$; therefore
(\ref{eq:L-g-K}) holds 
by the uniqueness property
of the bilateral Laplace transform. This implies the lemma statement.
\epr
\subsection*{Proof of Lemma~\ref{lem:LL}}
(a). For $s\in (0, \kappa_g^+)$, $a>0$ and $|\lambda|=1$ we have 
\begin{eqnarray*}
 \big[1-e^{-a(s+i\omega)}/\lambda\big]^{-m}
 = \sum_{i_1=0}^\infty\cdots \sum_{i_{m}=0}^\infty \big[e^{-a(s+i\omega)}/ \lambda\big]^{i_1+\cdots + i_m}
\\
 = \sum_{l=0}^\infty C_{l,m} \lambda^{-l} e^{-a l(s+i\omega)}.
\end{eqnarray*}
Therefore 
\begin{eqnarray*}
 \prod_{k=1}^q \big[1-e^{-a_k(s+i\omega)}/\lambda_k\big]^{-m_k}
 =\sum_{j_1=0}^\infty \cdots \sum_{j_q=0}^\infty 
 \prod_{k=1}^q C_{j_k, m_k}\lambda_k^{-j_k} e^{-a_kj_k(s+i\omega)},
\end{eqnarray*}
and it follows from  (\ref{eq:L-general}) that
\begin{align*}
 &L_{+,h}(t) = \frac{1}{2\pi}\int_{-\infty}^\infty  
 \frac{\widehat{K}((s+i\omega)h)\,\widehat{\psi}(-s-i\omega)}{\prod_{k=1}^q [1- e^{-a_k(s+i\omega)}/\lambda_k]^{m_k}}
 e^{(s+i\omega)t} \rd \omega
 \\
 &= \frac{1}{2\pi}\sum_{j_1=0}^\infty \cdots \sum_{j_q=0}^\infty 
 \prod_{k=1}^q C_{j_k, m_k}\lambda_k^{-j_k}
 \int_{-\infty}^\infty \widehat{K}((s+i\omega)h)\,\widehat{\psi}(-s-i\omega) 
 e^{
 (s+i\omega)(t-\sum_{k=1}^q a_kj_k)
 }
 \rd \omega\\
 &=
 \frac{1}{2\pi}\sum_{j_1=0}^\infty \cdots\sum_{j_q=0}^\infty 
 \prod_{k=1}^q C_{j_k, m_k}\lambda_k^{-j_k}
 \int_{-\infty}^\infty \widehat{K}(i\omega h)\,\widehat{\psi}(-i\omega) 
 e^{i\omega(t-\sum_{k=1}^q a_kj_k)}
 \rd \omega
 = \sum_{\ell\in \sL} \cC^+_\ell R_h(t-\ell), 
\end{align*}
where the third line follows from analyticity of the integrand, 
and $\cC_\ell^+$ is defined in (\ref{eq:C+1}).
Note that the change of the order of integration and summation is permissible
under the premise of the lemma.
\par 
(b). If $s\in (\kappa_g^{-},0)$ then 
 \begin{multline*}
 \big[1-e^{-a(s+i\omega)}/\lambda\big]^{-m} = 
 (-1)^m\lambda^m e^{am(s+i\omega)}[1- \lambda e^{a(s+i\omega)}]^{-m} 
 \\
 =(-1)^m
 \sum_{i_1=0}^\infty\cdots \sum_{i_{m}=0}^\infty \big[\lambda e^{a(s+i\omega)}\big]^{m+i_1+\cdots + i_m}
 = 
 (-1)^m\sum_{l=0}^\infty C_{l,m} \lambda^{l+m} e^{a (l+m)(s+i\omega)},
\end{multline*}
and, similarly to the above,
 \begin{eqnarray*}
 \prod_{k=1}^q \big[1-e^{-a_k(s+i\omega)}/\lambda_k\big]^{-m_k}
 =\sum_{j_1=0}^\infty \cdots \sum_{j_q=0}^\infty 
 \prod_{k=1}^q (-1)^{m_k}C_{j_k, m_k}\lambda_k^{j_k+m_k} e^{a_k(j_k+m_k)(s+i\omega)},
\end{eqnarray*}
which yields
\begin{eqnarray*}
 L_{s, h}(t)=L_{-,h}(t)
:=
 \sum_{j_1=0}^\infty \cdots \sum_{j_q=0}^\infty \Big[\prod_{k=1}^q (-1)^{m_k}C_{j_k, m_k} \lambda_k^{j_k+m_k}\Big] R_h\Big(t+ \sum_{k=1}^q a_k(j_k+m_k)\Big)
\\
=
\sum_{\ell\in \sL}^\infty \cC^-_\ell R_h\Big(t+\ell + \sum_{k=1}^qa_k m_k\Big),
 \end{eqnarray*}
 where $\cC_\ell^-$ is defined in (\ref{eq:C-1}).
 This completes the proof.
\epr
\subsection*{Proof of Theorem~\ref{th:upper}}
Throughout the proof we keep track of dependence of all constants on parameters of the 
classes $\sF^\prime_{\alpha, p}(A, B)$ and $\sG_{\alpha, p}^\prime(A,B)$. In what follows 
$c_1, c_2, \ldots$ stand for constants that can depend on parameters 
appearing in Assumptions~\ref{as:general}, \ref{as:psi} and ~\ref{as:C+-}, 
and on $\alpha$ and $p$ only.
For the sake of brevity, in the subsequent proof
we do not indicate integration limits if the corresponding integrals are taken over the entire real line.
\paragraph{Proof of statement~(a)} We assume that $x_0\geq 0$ and  consider the estimator 
$\tilde{f}_{+,h}^{(N)}(x_0)$ only; the derivation for $x_0<0$ and 
$\tilde{f}_{-,h}^{(N)}(x_0)$ is similar.
\par 
1$^0$. First we verify  
that under Assumption~\ref{as:psi} and condition~(K) 
the estimator $\tilde{f}_{+,h}^{(N)}(x_0)$
is well defined. 
 Because 
$K$ has finite support and it is infinitely differentiable on the real line, 
$\widehat{K}(i\omega)$ is also infinitely differentiable 
and rapidly decreasing as $|\omega|\to \infty$ in the sense that  
$|\omega|^k |\widehat{K}^{(j)}(i\omega)| \leq c(k,j)$ for all $k$ and $j$.
In particular, 
for  $k>\gamma+1$  by (\ref{eq:smooth}) of Assumption~\ref{as:psi}
we have 
\[
\int |\widehat{K}(i\omega h)|\, |\widehat{\psi}(-i\omega)| \rd \omega \leq  
c_2 h^{-k} <\infty.
\]
Thus function $R_h(\cdot)$ in (\ref{eq:Rh}) 
and kernels $L_{+,h}^{(N)}$
and $L_{-,h}^{(N)}$ in (\ref{eq:L+-})  
are well defined. 
\par 
2$^0$. First we derive an upper bound on the bias of $\tilde{f}_{+, h}^{(N)} (x_0)$. 
We have 
\begin{align}
 &\bE_f  \big[\tilde{f}_{+,h}^{(N)}(x_0)\big]  = \sum_{\ell\in\sL_N}\mathcal{C}_\ell^{+}
 \int R_{h}(y-x_{0}-\ell)f_{Y}(y)\,\rd y,
 \nonumber
 \\
  &\;= \sum_{\ell\in\sL_N}\mathcal{C}_\ell^{+}\cdot \frac{1}{2\pi}
 \int \widehat{K}(i\omega h)\widehat{\psi}(-i\omega)e^{-i\omega(x_{0}+\ell)}\left[
 \int  e^{i\omega y}f_{Y}(y)\,\rd y\right]\,\rd \omega 
\nonumber
 \\
  &\;= \frac{1}{2\pi}
 \int\widehat{K}(i\omega h)
 \widehat{\psi}(-i\omega)e^{-i\omega x_{0}}\widehat{f}_Y(-i\omega)
 \bigg[\sum_{\ell\in \sL_N} \cC_\ell^+ e^{-i\omega \ell}\bigg] \rd \omega~.
 \label{eq:E1}
 \end{align}
 It follows from definition of $\cC_\ell^+$ [cf.~(\ref{eq:C+1})] that 
 \begin{align*}
 &\sum_{\ell\in \sL_N} \cC_\ell^+ e^{-i\omega \ell} = \sum_{\ell \in \sL_N}
  \sum_{j_1=0}^\infty \cdots \sum_{j_q=0}^\infty {\bf 1}\{\ell = a^T j\} \bigg[\prod_{k=1}^q C_{j_k, m_k}
  \lambda_k^{-j_k}\bigg] e^{-i\omega \ell}
  \\
  &= \sum_{j_1=0}^N \cdots \sum_{j_q=0}^N 
   \bigg[\prod_{k=1}^q C_{j_k, m_k}
  (e^{-i \omega a_k}/\lambda_k)^{j_k}  \bigg] = \prod_{k=1}^q \sum_{j_k=0}^N C_{j_k, m_k}
  (e^{-i \omega a_k}/\lambda_k)^{j_k}. 
 \end{align*}
 Now noting that 
 \[
  \sum_{j=0}^N C_{j,m} [e^{-ia\omega }/\lambda]^j= \sum_{l_1=0}^N \cdots \sum_{l_m=0}^N 
  \bigg( \frac{e^{-i a\omega}}{\lambda}
  \bigg)^{l_1+\cdots+l_m} = \bigg[\frac{1- (e^{-i a \omega}/\lambda)^{N+1}}{1-(e^{-i a\omega}/\lambda)}\bigg]^m
 \]
we obtain 
\begin{align*}
\sum_{\ell\in \sL_N} \cC_\ell^+ e^{-i\omega \ell} = \prod_{k=1}^q 
 \bigg[\frac{1- (e^{-i a_k \omega}/\lambda_k)^{N+1}}{1-(e^{-i a_k\omega}/\lambda_k)}\bigg]^{m_k}.
\end{align*}
Substituting this expression in (\ref{eq:E1}), and taking into account (\ref{eq:factorization}) 
and  
$\widehat{f}_Y(-i\omega)=\widehat{f}(-i\omega) \widehat{g}(-i\omega)$ we obtain
\begin{align}
\bE_f  \big[\tilde{f}_{+,h}^{(N)}(x_0)\big]&  = \frac{1}{2\pi} \int
e^{-i\omega x_0} \widehat{K}(i\omega h) \widehat{f}(-i\omega) \prod_{k=1}^q \Big[1- 
\big(e^{-i a_k\omega}/\lambda_k\big)^{N+1}\Big]^{m_k} \rd \omega
\nonumber
\\
 &= \frac{1}{2\pi} \int e^{-i\omega x_0} \widehat{K}(i\omega h) \widehat{f}(-i\omega) 
\prod_{k=1}^q \sum_{j=0}^{m_k} \tbinom{m_k}{j} (-1)^{j} 
\big(e^{-i a_k\omega}/\lambda_k\big)^{j(N+1)}
\rd \omega 
\nonumber
\\
&= \sum_{j_1=0}^{m_1}\cdots \sum_{j_q=0}^{m_q} 
 \prod_{k=1}^q \bigg[\frac{ 
 (-1)^{j_k}\tbinom{m_k}{j_k}}{
 \lambda_k^{j_k (N+1)}}\bigg]
 \frac{1}{2\pi} \int
 \widehat{K}(i\omega h)\widehat{f}(-i\omega)
e^{- i\omega[x_0+(N+1)\sum_{k=1}^q j_k a_k]}\;
 \rd \omega 
\nonumber
\\
 &=\sum_{j_1=0}^{m_1}\cdots \sum_{j_q=0}^{m_q}
 \prod_{k=1}^q \bigg[\frac{ 
 (-1)^{j_k}\tbinom{m_k}{j_k}}{
 \lambda_k^{j_k (N+1)}}\bigg]
\frac{1}{h}\int 
 K\bigg(\frac{t-x_0-(N+1)\sum_{k=1}^q j_k a_k}{h}\bigg) f(t) \rd t
\nonumber
 \\ 
&= \frac{1}{h}\int K\Big(\frac{t-x_0}{h}\Big) f(t)\rd t + T_N(f;x_0),
\label{eq:E2}
 \end{align}
where we denote
\begin{align*}
 T_N(f;x_0) &:= \sum_{j_1=0}^{m_1}\cdots \sum_{j_q=0}^{m_q} 
 {\bf 1}\Big\{\sum_{k=1}^q j_k>0\Big\}
 \prod_{k=1}^q
\bigg[ \frac{
 (-1)^{j_k} \tbinom{m_k}{j_k}}
{ \lambda_k^{j_k (N+1)}}\bigg]\;
\\
 &\hspace{35mm}\times\;\frac{1}{h}\int 
 K\bigg(\frac{t-x_0-(N+1)\sum_{k=1}^q j_k a_k}{h}\bigg) f(t) \rd t.
\end{align*}
If $f\in \sM_{p} (B)$ then 
\begin{eqnarray}
 |T_N(f;x_0)| \leq c_1
 \sum_{j_1=0}^{m_1}\cdots \sum_{j_q=0}^{m_q}{\bf 1}\Big\{\sum_{k=1}^q j_k>0\Big\}
 \prod_{k=1}^q \tbinom{m_k}{j_k} 
\int_{-1}^1 f\Big(yh+x_0+(N+1)\sum_{k=1}^q j_ka_k\Big)\rd y 
\nonumber
\\
\leq \frac{c_2 B}{h(x_0+a_{\min}N)^{p}},
\label{eq:T-N}
\end{eqnarray}
where we denoted $a_{\min}:=\min\{a_1, \ldots, a_q\}$
and took into account that $a_{\min}>h$ for large enough $n$.
Then (\ref{eq:E2}), (\ref{eq:T-N}), condition~(K) and the fact that $k_0\geq \alpha +1$  
imply
the following upper bound on the bias
\begin{equation}\label{eq:bias-bound}
\sup_{f\in \sH_\alpha (A) \cap \sM_{p}(B)} \Big|
\bE_f[\tilde{f}_{+,h}^{(N)}(x_0)] - f(x_0)\Big| \leq c_3 
\bigg\{ Ah^\alpha + \frac{B}{h(x_0+a_{\min}N)^p}\bigg\}. 
\end{equation}
\par 
3$^0$. Now we bound the variance of $\tilde{f}_{+,h}^{(N)}(x_0)$.
We need the following notation.
For non--negative integer number $j$ 
we let 
\begin{equation}\label{eq:H-N-r}
H_{N, j} (\omega; x_0):= \left\{ \begin{array}{ll}
                             {\displaystyle \sum_{\ell\in \sL^*_N} 
                             \frac{\cC_\ell^+ e^{-i\omega(x_0+\ell)}}{(x_0+\ell)^{j}}
                             }, & x_0\geq 0,\\*[6mm]
                             {\displaystyle
                             \sum_{\ell\in \sL^*_N}  
                             \frac{\cC_\ell^- e^{-i\omega(x_0-\ell)}}{(x_0-\ell)^{j}},} & x_0<0,
                            \end{array}
\right.
\end{equation}
where 
$\sL^*_N:=\sL_N\backslash \{0\}$, and 
$\{\cC_\ell^+, \ell\in \sL\}$ and $\{\cC_\ell^{-}, \ell\in \sL\}$ are given in 
(\ref{eq:C+1}) and (\ref{eq:C-1}).
\par 
We have 
\begin{align}
  &{\rm var}_f\big\{\tilde{f}_{+,h}^{(N)}(x_0)\big\} \leq \frac{1}{n} \int 
 |L_{+,h}^{(N)}(y-x_0)|^2 f_Y(y) \rd y
\nonumber
 \\
 &\;\;= \frac{1}{n} \sum_{\ell_1\in \sL_N} \sum_{\ell_2\in \sL_N} 
 \cC^+_{\ell_1} \cC_{\ell_2}^+ \int R_h(y-x_0-\ell_1) \overline{R_h(y-x_0-\ell_2)}
 f_Y(y) \rd y
\nonumber
 \\
 &\;\;= \frac{1}{n} \int_{-\infty}^\infty |R_h(y-x_0)|^2 f_Y(y) \rd y
+\frac{2}{n} \sum_{\ell\in \sL_N^*} \cC_{\ell}^+ \int R_h(y-x_0-\ell) \overline{R_h(y-x_0)}
 f_Y(y) \rd y
\nonumber
 \\
&\;\;\;\hspace{15mm}+\frac{1}{n} \sum_{\ell_1\in \sL_N^*} \sum_{\ell_2\in \sL_N^*} 
 \cC^+_{\ell_1} \cC_{\ell_2}^+ \int R_h(y-x_0-\ell_1) \overline{R_h(y-x_0-\ell_2)}
 f_Y(y) \rd y 
\nonumber
 \\
 &\;\;=: \frac{1}{n} \int |R_h(y-x_0)|^2 f_Y(y) \rd y
 + \frac{2}{n} \sum_{\ell\in \sL_N^*} \cC_{\ell}^+ J_h(x_0+\ell, x_0)  
 \nonumber
 \\
&\;\;\;\;\hspace{10mm}+ 
 \frac{1}{n}\sum_{\ell_1\in \sL_N^*} \sum_{\ell_2\in \sL_N^*} 
 \cC^+_{\ell_1} \cC_{\ell_2}^+  J_h(x_0+\ell_1, x_0+\ell_2) =: S_1(x_0) + S_2(x_0) + S_3(x_0),
 \label{eq:var-decomp}
 \end{align}
where for $x_1, x_2>0$ we put
\begin{equation}\label{eq:J-h}
 J_h(x_1,x_2):=\frac{1}{4\pi^2}\iint  
 \widehat{R}_h(i\omega) \widehat{R}_h(-i\mu) 
\widehat{f}_Y(i(\omega-\mu)) e^{-i\omega x_1} e^{i\mu x_2} \rd\omega \rd\mu,
\end{equation}
where 
$\widehat{R}_h(i\omega)=\widehat{K}(i\omega h) \widehat{\psi}(-i\omega)$.
Since $\widehat{K}(i\omega)$ is rapidly decreasing as $|\omega|\to\infty$ and in view of 
(\ref{eq:smooth-1}) we have that $|\widehat{R}^{(j)}_h(i\omega)|\to 0$ as $|\omega|\to \infty$
for all $0\leq j\leq p_0$, where $p_0$ appears in Assumption~\ref{as:psi}.
Recall that by premise of the theorem $p_0\geq p$.
\par 
Now,  we proceed with bounding  the terms  $S_i(x_0)$, $i=1,2,3$ 
on the right hand side of (\ref{eq:var-decomp}).
First, since  $f\in \sM_p(B)$,  $f(x)\leq B$ and therefore 
$f_Y(y)\leq B$ for all $y$. By Parseval's identity and Assumption~\ref{as:psi}
\begin{align}
S_1(x_0) 
&=\frac{1}{n} \int |R_h(y-x_0)|^2 f_Y(y) \rd y 
\leq \frac{B}{4\pi^2n} 
\int|\widehat{R}_h(i\omega)|^2
 \rd \omega \leq  \frac{c_1B}{nh^{2\gamma+1}}.
 \label{eq:R-norm}
\end{align}
Furthermore, 
$f\in \sM_p(B)$ implies that $|\widehat{f}^{(j)}(i\omega|\leq B$, $0\leq j\leq p$, and since $\widehat{g}$
is analytic, 
the derivatives  $\widehat{f}_Y^{(j)}(i\omega)$ are finite for  $0\leq j\leq p$.
Therefore  for any integer $r\leq p$
by repeated integration by parts with respect to $\omega$ in (\ref{eq:J-h}) we obtain 
for $x_1\ne 0$
\begin{eqnarray*}
 J_h(x_1,x_2)= \frac{1}{4\pi^2} \int \widehat{R}_h(-i\mu) e^{i\mu x_2} 
 \int \frac{\rd^r}{\rd \omega^r} \Big[\widehat{R}_h(i\omega) \widehat{f}_Y (i(\omega-\mu)\Big] \frac{e^{-i\omega x_1}}{(-ix_1)^r} \rd \omega 
 \rd \mu
 \\
 = \frac{(-1)^r}{4\pi^2}\sum_{j=0}^r \tbinom{r}{j} \iint  
 \widehat{R}_h^{(j)}(i\omega) \widehat{f}_Y^{(r-j)} (i(\omega-\mu)) \widehat{R}_h(-i\mu)
 \frac{e^{-i\omega x_1}}{x_1^r} e^{i\mu x_2} \rd \omega \rd \mu.
\end{eqnarray*}
In the first line we have taken into account that $|\widehat{R}^{(j)}_h(i\omega)|\to 0$
as $|\omega|\to\infty$ for all $0\leq j\leq p_0$, and $p_0\geq p\geq r$. Now, invoking (\ref{eq:H-N-r})
we have
\begin{equation*}  
 S_2(x_0) 
= \frac{(-1)^r}{2\pi^2n}\sum_{j=0}^r \tbinom{r}{j} 
\iint  \widehat{R}_h^{(j)}(i\omega) 
\widehat{f}_Y^{(r-j)} (i(\omega-\mu)) \widehat{R}_h(-i\mu)
 H_{N,r}(\omega;x_0) e^{i\mu x_0} \rd \omega \rd \mu.
 \end{equation*}
The Cauchy--Schwarz inequality applied to the double integral on the right hand side yields
\begin{eqnarray}
&&S_2(x_0) \leq  \frac{1}{2\pi^2n}\sum_{j=0}^r  \tbinom{r}{j} 
\Big[\iint 
 \big|\widehat{R}_h^{(j)}(i\omega) H_{N,r}(\omega; x_0)\big|^2 \big|\widehat{f}_Y^{(r-j)}(i(\omega-\mu))\big| \rd\omega\rd \mu\Big]^{1/2}
\nonumber
 \\
&&\;\;\;\;\;\hspace{20mm}\times\; \Big[\iint 
\big|\widehat{R}_h(-i\mu)\big|^2  \big|\widehat{f}_Y^{(r-j)}(i(\omega-\mu))\big| \rd\omega\rd \mu\Big]^{1/2}
\label{eq:S22}
\\
&&\leq \frac{c_2}{n} \Big\{\max_{j=0,\ldots, r} \int 
|\widehat{f}_Y^{(r-j)}(i\omega)|\rd \omega \Big\}
\Big\{ \max_{j=0,\ldots, r} \int
\Big[|\widehat{R}_h^{(j)}(i\omega)|^2\big(1+ |H_{N,r}(\omega;x_0)|^2\big) \rd \omega \Big\}.
\nonumber
\end{eqnarray}
\par 
We bound $S_3(x_0)$ similarly. In particular, for any integer $r$ such that $2r\leq p$ 
repeated integration by parts in (\ref{eq:J-h}) first with respect to $\omega$ and then 
with respect to $\mu$ yields  
\begin{align*}
 J_h&(x_1, x_2) = \frac{1}{4\pi^2}
 \iint  
 \frac{\rd^r}{\rd \omega^r}\Big[ \widehat{R}_h(i\omega) \widehat{f}_Y(i(\omega-\mu))\Big]
 \widehat{R}_h(-i\mu) \frac{e^{-i\omega x_1}}{(-i)^r x_1^r} e^{i\mu x_2} \rd \omega \rd \mu
\\
& = \frac{(-1)^r}{4\pi^2} \sum_{j=0}^r \tbinom{r}{j} \int \widehat{R}^{(j)}_h(i\omega)
\frac{e^{-i\omega x_1}}{ x_1^r} \int 
\frac{\rd^r}{\rd \mu^r} 
\Big[ \widehat{f}_Y^{(r-j)}(i(\omega-\mu)) \widehat{R}_h(-i\mu)\Big] \frac{e^{i\mu x_2}}{(ix_2)^r} \rd \mu
\rd \omega
\\
&= \frac{1}{4\pi^2}
\sum_{j=0}^r \sum_{l=0}^r \tbinom{r}{j}\tbinom{r}{l}
\iint 
\widehat{R}_h^{(j)}(i\omega) \widehat{R}_h^{(l)}(-i\mu) 
\widehat{f}_Y^{(2r-j-l)} (i(\omega-\mu))
\frac{e^{-i\omega x_1} e^{i\mu x_2}}{x_1^r x_2^r} \rd \omega \rd \mu.
 \end{align*}
Hence 
\begin{align}
 &S_3(x_0)= \frac{1}{4\pi^2 n} \sum_{j=0}^r \sum_{l=0}^r \tbinom{r}{j}\tbinom{r}{l}
\nonumber
 \\
 &\;\times
 \iint 
\widehat{R}_h^{(j)}(i\omega) \widehat{R}_h^{(l)}(-i\mu) 
\widehat{f}_Y^{(2r-j-l)} (i(\omega-\mu))
H_{N,r}(\omega;x_0) H_{N,r}(-\mu;x_0) \rd \omega \rd \mu,
\label{eq:S3}
\end{align}
and by the Cauchy--Schwarz inequality 
\begin{align*}
 S_3(x_0)\;\leq\; &\frac{1}{4\pi^2 n} \sum_{j=0}^r \sum_{l=0}^r \tbinom{r}{j}\tbinom{r}{l}
\int 
\big|\widehat{f}_Y^{(2r-j-l)} (i\mu)\big| \rd \mu\; 
\\
&\;\;\times\;
\bigg(\int \big|\widehat{R}_h^{(j)}(i\omega)
H_{N,r}(\omega;x_0)\big|^2 \rd \omega\bigg)^{1/2}  
\bigg(\int \big|\widehat{R}_h^{(l)}(i\mu)
H_{N,r}(\mu;x_0)\big|^2 \rd \mu\bigg)^{1/2}
\\
&\leq   \frac{c_3}{n} 
\bigg\{\max_{j=0,\ldots, 2r} \int \big|\widehat{f}_Y^{(j)}(i\omega)\big|\rd \omega
\bigg\}
\bigg\{
\max_{j=0,\ldots, r} \int \big|\widehat{R}_h^{(j)}(i\omega) H_{N, r}(\omega;x_0)\big|^2 \rd \omega
\bigg\}. 
\end{align*}
Combining the above bounds on $S_1$, $S_2$ and $S_3$ we obtain that for any integer number $r$ such that 
$2r\leq p$ one has
\begin{eqnarray*}
 && {\rm var}_f\big\{\tilde{f}_{+,h}^{(N)}(x_0)\big\} 
 \leq \frac{c_1B}{nh^{2\gamma+1}} 
 \\
 &&\;\;\; +\; 
 \frac{c_4}{n} \Big\{\max_{j=0,\ldots, 2r} \int \big|\widehat{f}_Y^{(j)}(i\omega)\big|\rd \omega
\Big\}
\Big\{
\max_{j=0,\ldots, r} \int
\big|\widehat{R}_h^{(j)}(i\omega) \big|^2\big(1+ |H_{N, r}(\omega;x_0)\big|^2\big) \rd \omega
\Big\}.
\end{eqnarray*}
Now we bound the integrals on the right hand side of the above display formula.
\par 
Note that for any $j=0,\ldots, 2r,$ 
\[
 \int |\widehat{f}_Y^{(j)}(i\omega)|\rd \omega \leq c_5 \max_{l=0,\ldots, j} 
 \int |\widehat{g}^{(l)}(i\omega)|\cdot |\widehat{f}^{(j-l)}(i\omega)| \rd \omega.
\]
In view of  (\ref{eq:factorization}) for $l=0,\ldots, 2r$
\[
\widehat{g}^{(l)}(z) =\sum_{k=0}^l \tbinom{l}{k} 
\big[1/\widehat{\psi}(z)\big]^{(k)} \widehat{\varphi}^{(l-k)}(z),\;\;\;
\widehat{\varphi}(z):= \prod_{k=1}^q (1- e^{a_kz}\lambda_k^{-1})^{m_k}.
\]
It is obvious that $|\widehat{\varphi}^{(k)}(i\omega)|\leq c_6$ for all $k=0,\ldots, 2r$.
Furthermore, by the Fa\'a~di~Bruno formula
\[
  \big[1/\widehat{\psi}(z)\big]^{(k)}= 
  \sum_{m=1}^k \frac{(-1)^m m!}{[\widehat{\psi}(z)]^{m+1}}
  B_{k,m}\big(\widehat{\psi}^\prime(z), \ldots, \widehat{\psi}^{(k-m+1)}(z)\big),
\]
where $B_{k,m}(\cdot)$ are the Bell polynomials. Recall that $B_{k,m}$ is  a 
homogeneous polynomial in $k$ variables of degree $m$. 
This fact and (\ref{eq:smooth-1})  
imply that for  $k, m=0,\ldots, 2r,$ 
\[
 \big| B_{k,m}\big(\widehat{\psi}^\prime(i\omega), \ldots, \widehat{\psi}^{(k-m+1)}(i\omega)\big)\big| \leq 
 c_7 (1+|\omega|^\gamma)^m,\;\;\;\forall\omega\in \bR.
\]
Using  (\ref{eq:smooth}) we obtain  
$|[1/\widehat{\psi}(i\omega)]^{(k)}| \leq c_8 |\omega|^{-\gamma}$ for $|\omega|\geq \omega_0$; hence
$|\widehat{g}^{(l)}(i\omega)|\leq c_9 |\omega|^{-\gamma}$, $\forall|\omega|\geq \omega_0$.
This inequality in conjunction with boundedness of $\widehat{g}(i\omega)$ 
and $|\widehat{f}^{(j)}(i\omega)|\leq B$, $0\leq j\leq p$ 
for all $\omega$  implies that for
$l=0,\ldots, j$ and $j=0,\ldots 2r,$
\begin{eqnarray*}
 \int |\widehat{g}^{(l)}(i\omega)|\cdot |\widehat{f}^{(j-l)}(i\omega)| \rd \omega =
 \int_{|\omega|\leq \omega_0} + \int_{|\omega|\geq \omega_0}  |\widehat{g}^{(l)}(i\omega)|\cdot |\widehat{f}^{(j-l)}(i\omega)| \rd \omega
 \\
 \leq c_{10}B + c_{11} \int \frac{\widehat{f}^{(j-l)}(\omega)}{1+|\omega|^\gamma}\rd \omega. 
\end{eqnarray*}
Thus,  
if $f\in \sM^\prime_{p}(B)$ with $p\geq 2r$ then 
\[
 \max_{j=0, \ldots, 2r} 
 \int |\widehat{f}_Y^{(j)} (i\omega)|\rd \omega \leq c_9 \max_{j=0,\ldots, 2r} 
 \int \frac{|\widehat{f}^{(j)}(i\omega)|}{1+|\omega|^\gamma} \,\rd \omega  \leq c_{12} B,
\]
and therefore  
\[
 {\rm var}_f\big\{\tilde{f}_{+,h}^{(N)} (x_0)\big\} \leq 
 \frac{c_{13}B}{n} 
 \bigg\{
 \frac{1}{h^{2\gamma+1}} + 
  \max_{j=0,\ldots, r} \int \big|\widehat{R}_h^{(j)}(i\omega)\big|^2
  \big(1+ \big|H_{N, r}(\omega;x_0)\big|^2\big) \rd \omega
\bigg\}.
 \]
Taking into account that 
\[
 \widehat{R}_h^{(j)}(i\omega)=\frac{\rd^j}{\rd \omega^j} \Big[\widehat{K}(i\omega h) \widehat{\psi}(-i\omega)\Big]
 =\sum_{l=0}^j 
 \tbinom{j}{l} (ih)^l \widehat{K}^{(l)}(i\omega h) (-i)^{j-l}\widehat{\psi}^{(j-l)}(-i\omega)
\]
we have 
\begin{multline*}
 \int
 \big|\widehat{R}_h^{(j)}(i\omega)\big|^2\big(1+ |H_{N, r}(\omega;x_0)\big|^2\big) \rd \omega 
 \\
 \leq 
 c_{14}\sum_{l=0}^j h^l \int \big|\widehat{K}^{(l)}(i\omega h) \widehat{\psi}^{(j-l)}(-i\omega)\big|^2
 \big(1+|H_{N,r}(\omega;x_0)\big|^2\big)\rd \omega.
\end{multline*}
\par 
4$^0$. 
Combining the above bounds on the bias and variance we obtain that 
for any non--negative integer number $r$ satisfying $2r\leq p$ one has 
 \begin{multline}
  \cR_{n, \Delta_{x_0}}\big[\tilde{f}_{h}^{(N)}; \sF^\prime_{\alpha, p}(A, B)\big] 
\leq  c_{15} \Big\{  Ah^{\alpha} 
  + \frac{B}{h(x_0+a_{\min}N)^{p}} \Big\}
  \\*[2mm]
   \;\;\;\;+
  c_{16} \sqrt{\frac{B}{n}}\Big\{ \frac{1}{h^{2\gamma+1}}+ \sum_{l=0}^r h^{l} 
 \int \big| \widehat{K}^{(l)}(i\omega h)
\widehat{\psi}^{(r-l)}(-i\omega)\big|^2\big(1+ \big|H_{N, r}(\omega; x_0)\big|^2\big) \rd \omega \Big\}^{1/2},  
\label{eq:upper-bound}
\end{multline}
where
$C_1$ and $C_2$ may depend on $\alpha$ and $p$ only.
\par 
To complete the proof of statement~(a)
it is suffices to note that 
under Assumption~\ref{as:C+-}, $|H_{N, r}(\omega;x_0)|\leq c_{17}$, provided that 
$r\geq  \nu$. 
Therefore, in view of Assumption~\ref{as:psi} and condition (K) 
the last term on the right hand side of (\ref{eq:upper-bound})
is bounded above by $c_{18}Bh^{-2\gamma-1}$. Then the announced result follows by substitution  of the values of  $h$ and~$N$
in inequality~(\ref{eq:upper-bound}). 
\epr
\paragraph{Proof of statement~(b)}
The proof uses pointwise bounds derived in the proof of statement~(a). 
\par 
1$^0$. To derive the upper bound on the integrated squared bias consider equality~(\ref{eq:E2}). First we note the standard bound 
\cite[Section~1.2.3]{Tsybakov}:
\begin{equation}\label{eq:standard-bias}
 \sup_{f\in \sS_\alpha(A)} \int_{-\infty}^\infty  
 \bigg[\int_{-\infty}^\infty K\bigg(\frac{t-x_0}{h}\bigg) [f(t) - f(x_0)]\rd t \bigg]^2 \rd x_0 \leq c_1 A^2 h^{2\alpha}.
\end{equation}
Moreover, if $2p-1>0$ then by (\ref{eq:T-N})
\[
 \int_{0}^\infty |T_N(f;x_0)|^2 \rd x_0 \leq c_2 B^2 h^{-2} \int_0^\infty \frac{\rd x_0}{(x_0+a_{\min}N)^{2p}} \leq \frac{c_3B^2}{ 
 h^{2} (a_{\min}N)^{2p-1}}.
\]
The same upper bound holds for the integral of the squared bias of the estimator over $x_0\in (-\infty, 0]$. Thus,
\[
 \int_{-\infty}^\infty \Big|\bE_f[\tilde{f}_h(x_0)]-f(x_0)\Big|^2 \rd x_0 \leq c_4 \bigg(A^2h^{2\alpha} + \frac{B^2}{ 
 h^{2} (a_{\min}N)^{2p-1}}\bigg).
\]
\par 
2$^0$. Now consider the variance term. We use the variance decomposition given in (\ref{eq:var-decomp}). It follows from (\ref{eq:R-norm})
that
\[
 \int_{-\infty}^\infty S_1(x_0)\rd x_0 = \frac{1}{n} \int_{-\infty}^\infty |R_h(x)|^2 \rd x = \frac{1}{4\pi^2 n} \int_{-\infty}^\infty 
 |\widehat{K}(i\omega h)|^2|\widehat{\psi}(-i\omega)|^2
 \rd \omega \leq  \frac{c_1}{nh^{2\gamma+1}}.
\]
Furthermore, 
we note that for $r>1/2$
\begin{eqnarray}
 \int_0^\infty \big|H_{N, r}(\omega; x_0)\big|^2 \rd x_0  
 \leq \sum_{\ell_1\in \sL_N^*} \sum_{\ell_2\in \sL_N^*}  
 |\cC_{\ell_1}^+|\,|\cC_{\ell_2}^+| \int_0^\infty \frac{\rd x_0}{(x_0+\ell_1)^r (x_0+\ell_2)^{r}}
 \nonumber
 \\
 \leq c_2 \Big(\sum_{\ell\in \sL_N^*} |\cC_\ell^+| \ell^{-r+1/2}\Big)^2,
\label{eq:H2}
 \end{eqnarray}
 and the same inequality holds for 
 the integral 
 $\int_{-\infty}^0 |H_{N, r}(\omega; x_0)|^2\rd x_0$.
By Assumption~\ref{as:C+-}, if $r\geq \nu+1/2$ then the sum on the right hand side of (\ref{eq:H2})
is uniformly bounded in $N$. 
Therefore
using (\ref{eq:S22}) and applying the Cauchy--Schwarz inequality we obtain
\begin{align*}
 \int_{-\infty}^\infty S_2(x_0)\rd x_0 \;
 &\leq\;  \frac{c_3}{n}  \Big(\sum_{\ell\in \sL_N^*} \big(|\cC_\ell^+| \vee |\cC_\ell^-|\big)\ell^{-r+1/2}\Big)^2
 \\
 & \;\;\;\times 
 \sum_{j=0}^r \tbinom{r}{j} \iint
 \Big|
 \widehat{R}_h^{(j)}(i\omega) \widehat{f}_Y^{(r-j)}(i(\omega-\mu)) \widehat{R}_h(-i\mu)  
 \Big|
 \rd \omega \rd \mu 
 \\
 &\leq  \frac{c_4}{n}  
 \Big\{\max_{j=0,\ldots, r} \int |\widehat{f}_Y^{(j)}(i\omega)|\rd \omega \Big\}
\Big\{ \max_{j=0,\ldots, r} \int \big|\widehat{R}_h^{(j)}(i\omega)\big|^2 \rd \omega \Big\}.
\end{align*}
The term originating from $S_3(x_0)$ is bounded similarly. 
The bound (\ref{eq:H2}) together with 
with (\ref{eq:S3})
yields 
\begin{align*}
 \int_{-\infty}^\infty S_3(x_0) \rd x_0 \;\leq  \;& \frac{c_5}{n}
\Big\{\max_{j=0,\ldots, 2r} \int \big|\widehat{f}_Y^{(j)}(i\omega)\big|\rd \omega
\Big\}
\Big\{
\max_{j=0,\ldots, r} \int \big|\widehat{R}_h^{(j)}(i\omega)\big|^2 \rd \omega
\Big\}. 
\end{align*}
\par 
3$^0$.
Combining the obtained inequalities with the bound on the bias and using the same reasoning as in the proof of Theorem~\ref{th:upper} 
we conclude that for $r\geq \nu+1/2$ we have 
\begin{eqnarray*}
 \cR_{n, \Delta_2}[\tilde{f}_h^{(N)}; \sG^\prime_{\alpha, p}(A, B)]\leq  c
 \Big\{A^2h^{2\alpha} + \frac{B^2}{ 
 h^{2} (a_{\min}N)^{2p-1}}  + \frac{B}{nh^{2\gamma+1}}\Big\}^{1/2}.
\end{eqnarray*}
Substitution of the values $h=h_*$ and $N$ completes the proof of statement~(b).

\epr

 \subsection*{Proof of Theorem~\ref{th:conv-uniform}}
 
In the subsequent proof we keep track of all constants depending on parameters of classes 
$\sF_{\alpha, p}(A, B)$ and $\sG_{\alpha, p}(A, B)$. In what follows  
$c_1, c_2, \ldots$ denote positive constants that depend on $m$ and $\alpha$ only.  
\paragraph{Proof of statement (a)}   We provide the proof  of statement (a) for  the estimator corresponding to $x_0\geq 0$ only. The proof for $x_0<0$ is identical in every detail.
The estimator is given by the formula 
 \[
 \tilde{f}_{+, h}^{(N)}(x_0)
 =\frac{1}{n}\sum_{k=1}^n \frac{(2\theta)^m}{h^{m+1}} \sum_{j=0}^N C_{j,m} K^{(m)} 
 \Big(\frac{Y_k-x_0 -\theta(2j+m)}{h}\Big).
\]
\par 
1$^0$. 
The variance of this estimator is bounded as follows:
\begin{align*}
 &{\rm var}_f \big[\tilde{f}_{+,h}^{(N)}(x_0)\big] \leq \frac{1}{n} \int_{-\infty}^\infty
 |L^{(N)}_{+,h}(y-x_0)|^2 f_Y(y) \rd y
 \\
 &=  \frac{(2\theta)^{2m}}{nh^{2m+2}}  
 \int_{-\infty}^\infty \Big|\sum_{j=0}^N C_{j,m}K^{(m)} \Big(
 \frac{y-x_0-\theta(2j+m)}{h}\Big)\Big|^2 f_Y(y) \rd y
\\
&= 
 \frac{(2\theta)^{2m}}{nh^{2m+2}} \sum_{j=0}^N 
\sum_{l=0}^N C_{j,m} C_{l,m}
\\
&\;\;\;\;\times\;
\int_{-\infty}^\infty 
K^{(m)} \Big(
 \frac{y-x_0-\theta(2j+m)}{h}\Big)K^{(m)} \Big(
 \frac{y-x_0-\theta(2l+m)}{h}\Big)
f_Y(y) \rd y.
 \end{align*}
 Assume that $h < \theta$ (this is always fulfilled for large $n$), and
denote 
\[
 I_j(x_0):=[x_0+\theta(2j+m)-h, x_0+\theta(2j+m)+h],\;\;\;j=0,1,2,\ldots.
\]
Since 
 ${\rm supp}(K)\subseteq [-1,1]$
and $h<\theta$, the intervals  
$I_j(x_0)$ and $I_l(x_0)$ are disjoint for $j\ne l$.
Therefore
\begin{eqnarray}
 {\rm var}_f\big[\tilde{f}_{+,h}^{(N)}(x_0)\big] &\leq& \frac{(2\theta)^{2m}}{nh^{2m+2}}
 \sum_{j=0}^N C_{j,m}^2 \int_{-\infty}^\infty 
 \Big|K^{(m)} \Big(
 \frac{y-x_0-\theta(2j+m)}{h}\Big)\Big|^2 f_Y(y) \rd y
 \nonumber
 \\
 &=& \frac{(2\theta)^{2m}}{nh^{2m+1}}  \sum_{j=0}^N C_{j,m}^2 
 \int_{-1}^1 [K^{(m)}(y)]^2 
 f_Y(x_0+\theta (2j+m)+yh) \rd y
\nonumber 
 \\
 &\leq&\frac{c_1 \theta^{2m}}{nh^{2m+1}}
  \sum_{j=0}^N \frac{C_{j,m}^2}{h} 
 \int_{I_j(x_0)} f_Y(t) \rd t,
\label{eq:var-conv-uniform}
 \end{eqnarray}
 where we have used that $K^{(m)}(\cdot)$ is bounded above by a constant.
Now we bound the sum on the right hand side of (\ref{eq:var-conv-uniform}).
\par 
First, 
note that $g$ is supported on $[-m\theta, m\theta]$, and $g(x)\leq c_2/\theta$. Therefore
\[
 f_Y(t)= \int_{-m\theta}^{m\theta} f(t-y) g(y)\rd y \leq \frac{c_2}{\theta} \int_{-m\theta}^{m\theta} f(t-y)\rd y. 
\]
Furthermore, writing $\xi_j:=x_0+(2j+m)\theta$ for brevity we have 
\begin{align*}
 &\frac{1}{h} \int_{I_j(x_0)} f_Y(y) \rd y \leq \frac{1}{h}\int_{-h}^h \frac{c_2}{\theta}\int_{-m\theta}^{m\theta}
 f(y+\xi_j-z)\rd z\rd y
 \\
 & = \frac{c_2}{\theta h} \int_{-\infty}^\infty \int_{-\infty}^\infty  f(y-z) {\bf 1}\{-m\theta \leq z\leq m\theta\}
 {\bf 1}\{\xi_j-h\leq y \leq \xi_j+h\}\rd z\, \rd y
 \\
 & = 
 \frac{c_2}{\theta h} \int_{-\infty}^\infty   f(u) \int_{-\infty}^\infty {\bf 1}\{u-m\theta \leq y\leq u+m\theta\}
 {\bf 1}\{\xi_j-h\leq y \leq \xi_j+h\}\rd y\, \rd u,
\end{align*}
and since $m\theta>h$ we obtain
\begin{align}
 \frac{1}{h} \int_{I_j(x_0)} f_Y(y) \rd y &= \frac{c_2}{\theta} \int_{-h}^h f(t+\xi_j+m\theta)
 \Big(1-\frac{t}{h}\Big)\rd t 
 \nonumber
 \\
 &\;\;+ \frac{c_2}{\theta} \int_{-h}^h f(t+\xi_j-m\theta)
 \Big(1+\frac{t}{h}\Big)\rd t
  + \frac{2c_2}{\theta} \int_{h-m\theta}^{-h+m\theta} f(t+\xi_j) \rd t
\nonumber
  \\
 &\leq\;\; 
 \frac{c_3}{\theta} \int_{-h}^h f(t+x_0+2(j+m)\theta)\rd t 
 + \frac{c_3}{\theta} \int_{-h}^h f(t+x_0+2j\theta)\rd t 
 \nonumber
 \\
 & \hspace{40mm} + \frac{c_3}{\theta} 
 \int_{-m\theta}^{m\theta} f(t+x_0+(2j+m)\theta) \rd t.
 \label{eq:integrals}
\end{align}
Note that
$C_{j,m}\leq \big[(j+m-1)/(m-1)\big]^{m-1} e^{m-1} \leq c_4 j^{m-1}$ for $j\geq 1$. 
Taking into account that
$f\in \sM_{2m-2}(B)$ we have 
\begin{align*}
 \sum_{j=0}^N  \frac{C_{j,m}^2}{\theta}\int_{-h}^h &f(t+x_0+2(j+m)\theta)\rd t
\\
 &
\leq 
 c_5 \sum_{j=0}^N \frac{j^{2m-2}}{\theta(x_0+ 2\theta j)^{2m-2}} \int_{x_0+2(j+m)\theta-h}^{x_0+2(j+m)\theta+h} 
t^{2m-2}f(t)\rd t
 \\
 &\leq 
\frac{c_6}{\theta^{2m-1}} \int_{x_0+2m\theta-h}^{x_0+2(N+m)\theta+h} t^{2m-2} f(t)\rd t\leq 
 c_6B\theta^{-2m+1},
\end{align*}
where 
the first line follows from  elementary inequality $\int_a^b f(t)\rd t\leq \int_a^b (t/a)^k f(t)\rd t$
for $a>0$, $k\geq 0$ and non--negative $f$, while  the second line   
follows from the fact that the integrals under the sum 
are taken over disjoint intervals because $h<\theta$.
The similar upper bound holds 
for the sum corresponding to the second integral on the right hand side of (\ref{eq:integrals}). 
The expression corresponding  to the  third integral is bounded as follows
\begin{eqnarray*}
 \sum_{j=0}^N  \frac{C_{j,m}^2}{\theta} \int_{x_0+2\theta j}^{x_0+(2j+2m)\theta} f(t) \rd t\leq 
 \sum_{j=0}^N  \frac{C_{j,m}^2}{\theta(x_0+2\theta j)^{2m-2}} \int_{x_0+2\theta j}^{x_0+(2j+2m)\theta} 
 t^{2m-2}f(t)  \rd t
 \\
\leq 
 \frac{c_7}{\theta^{2m-1}} \sum_{j=0}^N    \int_{x_0+2\theta j}^{x_0+(2j+2m)\theta} t^{2m-2} f(t) \rd t
 \leq c_8 B\theta^{-2m+1},
\end{eqnarray*}
where the last inequality holds because 
\[
\sum_{j=0}^N    \int_{x_0+2\theta j}^{x_0+(2j+2m)\theta} t^{2m-2} f(t) \rd t \leq m \int_{x_0}^{x_0+(2N+2m)\theta} t^{2m-2} f(t) \rd t \leq mB.
\]
Therefore 
\[
 {\rm var}_f \big[ \tilde{f}_{+,h}^{(N)}(x_0)\big] \leq c_9 B \theta (nh^{2m+1})^{-1}~. 
\]
\par 
2$^0$. Following   the proof of Theorem~\ref{th:upper} preceding formula (\ref{eq:T-N}) we have that
\[
 \bE_f \big[ \tilde{f}_{+,h}(x_0)\big]= \frac{1}{h} \int_{-\infty}^\infty 
 K\Big(\frac{t-x_0}{h}\Big) f(t)\rd t + T_N(f;x_0),
\]
where 
\begin{eqnarray}
     T_N(f;x_0)&=& \sum_{j=1}^m \tbinom{m}{j} (-1)^j \frac{1}{h}\int_{-\infty}^\infty K\Big(\frac{t-x_0-2\theta (N+1)j}{h}\Big)
     f(t)\rd t
     \nonumber
     \\
     &=&\sum_{j=1}^m \tbinom{m}{j} (-1)^j \int_{-1}^{1} K(y) f(yh+x_0+2\theta (N+1)j)\rd y.
   \label{eq:TN-conv}
   \end{eqnarray}
Letting $\xi_j= x_0+2\theta(N+1)j$ and taking into account that 
$f\in \sM_{p}(B)$,  and $\theta>h$ 
we have for any $j=1,\ldots, m$
\begin{equation}\label{eq:TN-conv1}
 \int_{-1}^{1} |K(y)| f(yh+x_0+2\theta (N+1)j)\rd y \leq \frac{c_{10}}{h} 
 \int_{\xi_j-h}^{\xi_j+h} \frac{t^{p}f(t)}{(\xi_j-h)^{p}}\rd t
\leq \frac{c_{10} B}{h (\theta N)^{p}}.
 \end{equation}
These inequality yields
$|T_N(f;x_0)| 
 \leq c_{11} B h^{-1} (\theta N)^{-p}$,
and 
since $f\in \sH_\alpha(A)$, 
\[
 \Big|\bE_f\big[\tilde{f}_{+,h}^{(N)}(x_0)\big] -f(x_0)\Big| \leq c_{12} Ah^\alpha + c_{11} B h^{-1} (\theta N)^{-p}.
\]
Then statement~(a) follows immediately from the established bounds 
on the bias and variance 
by substitution of the  chosen values of $h$ and $N$.
\paragraph{Proof of statement~(b)}
We start with the bounding the  variance term.
The basis for the derivation is formula 
(\ref{eq:var-conv-uniform}) that should be integrated over \mbox{$x_0\in [0,\infty)$}.
In view of (\ref{eq:integrals}) and the subsequent formulas
\begin{multline*}
\sum_{j=0}^N  \frac{C^2_{j,m}}{h} \int_{I_j(x_0)} f_Y(y) \rd y 
\\
\;\;\;\;\;\leq\; 
\sum_{j=0}^N \frac{c_1C_{j,m}^2}{\theta}
\bigg[
\int_{x_0+2\theta(m+j)-h}^{x_0+2\theta(m+j)+h} f(t)\rd t +
\int_{x_0+2\theta j-h}^{x_0+2\theta j+h} f(t)\rd t
+
\int_{x_0+2\theta j}^{x_0+2\theta (m+j)} f(t)\rd t\bigg]
\\
 =: J_1(x_0)+J_2(x_0)+J_3(x_0).
\end{multline*}
Since $f\in \sM_p(B)$ we have  
\begin{eqnarray*}
 J_1(x_0)& =& \sum_{j=0}^N \frac{c_1 C_{j,m}^2}{\theta (x_0+2\theta (m+j)-h)^p}
 \int_{x_0+2\theta(m+j)-h}^{x_0+2\theta (m+j)+h} t^pf(t)\rd t
\\
 &\leq& 
 \frac{c_1}{\theta (x_0+\theta)^{p-2m+2}}
 \sum_{j=0}^N \frac{C_{j,m}^2}{(x_0+2\theta j+ \theta (2m-1))^{2m-2}}
 \int_{x_0+2\theta(m+j)-h}^{x_0+2\theta (m+j)+h} t^pf(t)\rd t
 \\
 &\leq& 
 \frac{c_2}{\theta^{2m-1} (x_0+\theta)^{p-2m+2}}
 \sum_{j=0}^N 
 \int_{x_0+2\theta(m+j)-h}^{x_0+2\theta (m+j)+h} t^pf(t)\rd t\leq 
 \frac{c_3B}{\theta^{2m-1} (x_0+\theta)^{p-2m+2}}.
\end{eqnarray*}
Hence 
taking into account that $p>2m-1$ and integrating with respect to $x_0\in [0,\infty)$ we obtain 
\[
 \int_0^\infty J_1(x_0)\rd x_0 \leq c_3B\theta^{-p}.  
\]
Using the same reasoning it is immediate to show that the same bound holds for the integrals 
of $J_2(x_0)$ and $J_3(x_0)$:
$\int_0^\infty J_i(x_0)\rd x_0 \leq c_4B\theta^{-p}$, $i=1,2$.
Combining these bounds with (\ref{eq:var-conv-uniform}) we obtain
\[
 \int_0^{\infty} {\rm var}_f\big\{ \tilde{f}_h^{(N)}(x_0)\big\} \rd x_0 \leq  c_4 B\theta^{2m-p} (nh^{2m+1})^{-1}. 
\]
 The integral over the negative semi--axis is bounded similarly.
\par 
To bound the integrated squared bias we note that  (\ref{eq:TN-conv})--(\ref{eq:TN-conv1}) and $f\in \sM_{p}(B)$ 
imply
\[
 \int_0^\infty T_N^2(f; x_0)\rd x_0 \leq  c_1 B^2 h^{-2} (\theta N)^{-2p+1}, 
\]
and the same estimate holds for the integral of $T_N(f;x_0)$ over the negative semi--axis. In view of  
 (\ref{eq:standard-bias})
we obtain 
 \[
 \int_{-\infty}^\infty \Big|\bE_f[\tilde{f}_h(x_0)]-f(x_0)\Big|^2 \rd x_0 \leq c_2 \bigg(A^2h^{2\alpha} + \frac{B^2}{ 
 h^{2} (\theta N)^{2p-1}}\bigg).
\]
Then the statement follows from the established upper bounds on the integrated variance and the integrated squared bias. 
 \epr 
 
 \subsection*{Proof of Theorem~\ref{th:lower-bound-2}}
 
 In the subsequent proof $c_1, c_2, \ldots$ stand for positive constants that do not depend on $n$. 
The proof of (\ref{eq:low-b-1}) is based on the 
standard reduction to a two--point testing problem, while in the 
proof of (\ref{eq:low-b-2})
we use reduction to the problem of testing  
multiple hypotheses [see \cite[Chapter~2]{Tsybakov}]. 
 \paragraph{Proof of statement~(a)} 
Let $r>1/2$ be a real  number and consider the probability density
\begin{equation}\label{eq:f-0}
 f_0(x) := \frac{C_r}{(1+x^2)^r},\;\;\;x\in \bR,
\end{equation}
where $C_r$ is a normalizing constant. Clearly, 
$f_0\in \sM_p(B)$ for  $p<2r-1$ and sufficiently large constant $B$ depending on $p$. 
In addition, $f_0$ is infinitely differentiable and belongs to $\sH_\alpha(A)$
for any $\alpha$ and large enough $A$.
\par\medskip
Pick function $\widehat{\psi}_0$ with the following properties:
\begin{itemize}
\item[(i)] $\widehat{\psi}_0(\omega)=\widehat{\psi}_0(-\omega)$, $\forall \omega$;
\item[(ii)] $\widehat{\psi}_0(\omega)=1$ for $\omega\in [-1+\delta, 1-\delta]$ with some fixed $\delta>0$, 
and $\widehat{\psi}_0(\omega)=0$, $|\omega|>1$;
 \item[(iii)]  $\widehat{\psi}_0$  monotonically climbs from $0$ to $1$  
 on $[-1, -1+\delta]$. In addition,  $\widehat{\psi}_0$ is infinitely  differentiable function on the real line.
\end{itemize}
\par
Let $h>0$ be a small real number such that $h<\pi/\theta$, and $N\geq 1$ be an integer number.  Define
\[
 \widehat{\psi}_h(\omega):= \sum_{k=N+1}^{2N}
 \bigg[\widehat{\psi}_0\Big(\frac{\omega-\pi k/\theta}{h}\Big)+
 \widehat{\psi}_0\Big(\frac{\omega+\pi k/\theta}{h}\Big)\bigg].
\]
Note that $\widehat{\psi}_h$ is even 
and supported on 
the union of disjoint sets $\cup_{k=N+1}^{2N} A_k$, where 
$A_k:=[-\pi k/\theta-h, -\pi k/\theta+h]\cup [\pi k/\theta-h, \pi k/\theta+h]$.
Function $\psi_h$ is given by the inverse Fourier transform:
\[
 \psi_h(x)= \frac{1}{2\pi}\int_{-\infty}^\infty 
 \widehat{\psi}_h(\omega) e^{i\omega x}\rd \omega =
 2h\psi_0(hx) \sum_{k=N+1}^{2N}
 \cos\Big(\frac{\pi k x}{\theta}\Big),\;\;\;x\in \bR.
\]
For real numbers $M>0$ and $c_0>0$ define
\begin{equation}\label{eq:f-1}
 f_1(x):=f_0(x) + c_0 M \psi_h(x). 
\end{equation}
 \par 
1$^0$. We demonstrate that under appropriate choice of constants $M$, $h$ and $N$ function $f_1$ 
is a probability density, and 
it belongs to $\sH_\alpha(A)\cap \sM_p(B)$ with $p<2r-1$.
\par 
First, we note that $\int_{-\infty}^\infty \psi_h(x) \rd x =0$ because $\widehat{\psi}_h(0)=0$.
Second, since $\widehat{\psi}_0$ is 
infinitely differentiable, $\psi_0$ is rapidly decreasing, and 
$|\psi_0(x)|\leq c_1 |x|^{-2r}$  for some constant $c_1=c_1(r)$ depending on $r$.
Therefore $|\psi_h(x)| \leq c_2 |x|^{-2r} h^{-2r+1}N$ for all $x$, and 
if we set  
\begin{equation}\label{eq:M-h}
 M= h^{2r-1} N^{-1},
\end{equation}
then by choice of constant $c_0$ we can ensure that $c_0 M |\psi_h(x)|\leq f_0(x)$ for all $x$.
Thus, $f_1$ is a probability density provided  (\ref{eq:M-h}) holds. This also shows that $f_1\in \sM_p(B)$ for $p<2r-1$.
\par 
For simplicity assume that $\alpha$ is integer; then
\begin{eqnarray*}
 |\psi_h^{(\alpha)}(x)| = \bigg|2h \sum_{j=0}^\alpha \tbinom{\alpha}{j} h^j \psi_0^{(j)}(xh) 
 \sum_{k=N+1}^{2N}[\cos 
 \big(\pi k x/\theta\big)]^{(\alpha-j)}\bigg| 
 \\
 \leq c_2 h
 \sum_{j=0}^\alpha h^j N^{\alpha-j+1} \leq c_3  hN^{\alpha+1}.
\end{eqnarray*}
This implies that $f_1\in \sH_\alpha(A)$ if $MhN^{\alpha+1}=h^{2r} N^\alpha\leq A$.
\par 
2$^0$. Without loss of generality we consider the problem of estimating  the value $f(0)$. Note that we have 
\[
 |f_1(0)-f_0(0)|= c_0 M\psi_h(0)=c_0 \psi_0(0) MhN = c_4 h^{2r}.
\]
 \par 
3$^0$. Now we bound from above the $\chi^2$--divergence between densities of observations $f_{Y,0}$ and $f_{Y,1}$ corresponding to the hypotheses
$f=f_0$ and $f=f_1$. We have 
\begin{eqnarray*}
 \chi^2(f_{Y,0}; f_{Y,1}) = c_0^2M^2\int_{-\infty}^\infty \frac{|(g\ast \psi_h)(x)|^2}{(g\ast f_0)(x)} \rd x
\end{eqnarray*}
Since $g$ is supported on $[-m\theta, m\theta]$ we have 
\[
 (g\ast f_0)(x) = \int_{-m\theta}^{m\theta}  \frac{C_r g(y)}{ [1+(x-y)^2]^r}\rd y\geq \frac{c_5}{[1+(|x|-m\theta)^2]^r},\;\;\;\forall |x|>m\theta,
\]
and 
\[
(g\ast f_0)(x)\geq c_6, \;\;\;|x|\leq m\theta.     
\]
Therefore 
\begin{eqnarray*}
 \chi^2(f_{Y,0}; f_{Y,1}) \leq c_7 M^2\int_{-\infty}^\infty  |(g\ast \psi_h)(x)|^2 \rd x + 
 c_8 M^2 \int_{-\infty}^\infty   x^{2r}|(g\ast \psi_h)(x)|^2 \rd x
 \\
 =: c_7I_1 + c_8I_2.
\end{eqnarray*}
Now we bound integrals $I_1$ and $I_2$ on the right hand side. First we assume that $r$ is an integer number.
Then by Parseval's identity  
\begin{align*}
 &I_1 = \frac{M^2}{2\pi} \int_{-\infty}^\infty |\widehat{g}(\omega)|^2 |\widehat{\psi}_h(\omega)|^2\rd \omega 
 = \frac{M^2}{\pi} 
 \int_{-\infty}^{\infty} \bigg|\frac{\sin (\theta \omega)}{\theta \omega}\bigg|^{2m}\; 
\sum_{k=N+1}^{2N}
 \bigg|\widehat{\psi}_0\Big(\frac{\omega-\pi k/\theta}{h}\Big)\bigg|^2\rd \omega 
 \\
 &=
 \frac{M^2h}{\pi} \sum_{k=N+1}^{2N}
 \int_{-1}^{1} \bigg|\frac{\sin (\pi k+ \theta\xi h )}{\pi k+\theta \xi h}\bigg|^{2m}\; 
 \big|\widehat{\psi}_0(\xi)\big|^2\rd \xi  \leq 
 c_9 M^2 h^{2m+1} N^{-2m+1}=\frac{c_9 h^{2m+4r-1}}{N^{2m+1}}~.
\end{align*}
Furthermore,
\begin{equation}\label{eq:I22}
 I_2 =  M^2 \int_{-\infty}^\infty x^{2r}
 |(g\ast \psi_h)(x)|^2\rd x= 
 \frac{M^2}{2\pi}\int_{-\infty}^\infty
 \Big|\frac{\rd^r}{\rd \omega^r} \widehat{g}(\omega)\widehat{\psi}_h(\omega)\Big|^2\rd \omega.
 \end{equation}
 We have 
 \begin{eqnarray}
 && \frac{\rd^r}{\rd \omega^r} \widehat{g}(\omega)\widehat{\psi}_h(\omega)
 = \sum_{j=0}^r \tbinom{r}{j} \widehat{g}^{(j)}(\omega) \widehat{\psi}_h^{(r-j)}(\omega)
\nonumber
\\
 &&\;\;\;\;\;\;\;= \sum_{j=0}^r \tbinom{r}{j}  \frac{\widehat{g}^{(j)}(\omega)}{h^{r-j}} 
 \sum_{k=N+1}^{2N} 
 \bigg[\widehat{\psi}_0^{(r-j)}\Big(\frac{\omega- \pi k/\theta}{h}\Big) 
+ 
 \widehat{\psi}_0^{(r-j)}\Big(\frac{\omega+ \pi k/\theta}{h}\Big)\bigg].
\label{eq:dromega}
 \end{eqnarray}
 By the Fa\'a~di~Bruno formula for $j\geq 1$ 
 \begin{align*}
 \widehat{g}^{(j)}(\omega)&= 
 \frac{\rd^j}{\rd \omega^j} \bigg(\frac{\sin \theta \omega}{\theta \omega}\bigg)^m 
 \\
 &=\sum_{l=1}^j j(j-1)\cdots (j-l+1) \bigg(\frac{\sin \theta \omega}{\theta \omega}\bigg)^{m-l} B_{j,l}\Big(\widehat{g}_0^\prime(\omega),\ldots,
 \widehat{g}_0^{(j-l+1)}(\omega)\Big),
 \end{align*}
 where 
 $B_{j,l}(\cdot)$ are the Bell polynomials, 
 and $\widehat{g}_0(\omega):=(\sin \theta\omega)/(\theta\omega)$.
 First, we note that 
$|\widehat{g}_0^{(k)}(\omega)| \leq c_{10}(k) \min \{ |\omega|^{-1}, 1\}$ for any $k$. 
Using this fact and taking into account that  $B_{j, l}$ is a homogeneous polynomial in $j$ variables of degree $l$
we have 
\[
 |\widehat{g}^{(j)}(\omega)| \leq  c_{11} \sum_{l=1}^j   \bigg|\frac{\sin \theta \omega}{\theta \omega}\bigg|^{m-l} |\omega|^{-l} = \frac{c_{11}}{|\theta \omega|^m} 
 \sum_{l=1}^j |\sin \theta\omega|^{m-l}~.
\]
Taking into account this inequality, (\ref{eq:dromega}), (\ref{eq:I22}),  
the fact that $\widehat{\psi}_0$ is supported on  $[-1,1]$, 
and recalling that $A_k:=[-\pi k/\theta-h, -\pi k/\theta+h]\cup [\pi k/\theta-h, \pi k/\theta+h]$
 are disjoint for all $k=N+1,\ldots, 2N$, we obtain 
\begin{align}
 I_2 \;\leq\; 
 c_{12}M^2
 \int_{-\infty}^\infty \bigg|\sum_{j=0}^r \tbinom{r}{j} \frac{\widehat{g}^{(j)}(\omega)}{h^{r-j}}
 \sum_{k=N+1}^{2N} \Big\{\widehat{\psi}_0^{(r-j)}
 \Big(\frac{\omega-\pi k/\theta}{h}\Big)+\widehat{\psi}_0^{(r-j)}
 \Big(\frac{\omega-\pi k/\theta}{h}\Big)\Big\}\bigg|^2 \rd \omega
\nonumber
 \\
 \leq \; c_{13}M^2 h^{-2r}\sum_{k=N+1}^{2N} \int_{A_k} \Big|\sum_{j=0}^r
 h^j\widehat{g}^{(j)}(\omega)\Big|^2 \rd \omega
\nonumber
 \\
 \leq\; c_{14} M^2h^{-2r} \sum_{k=N+1}^{2N} \int_{A_k} 
 \bigg(\Big|\frac{\sin \theta \omega}{\theta\omega}\Big|^{2m}+ 
 \frac{1}{|\theta\omega|^{2m}} \sum_{j=1}^r h^{2j} \sum_{l=1}^j |\sin\theta\omega|^{2m-2l}
 \bigg)\rd \omega
 \nonumber
 \\
  \leq c_{15} M^2 h^{-2r+1+2m} \sum_{k=N+1}^{2N} \frac{1}{k^{2m}} = c_{16} h^{2m+2r-1}N^{-2m-1}. 
  \label{eq:I222}
\end{align}
Combining these bounds we obtain that for integer $r\geq 1$
\begin{equation}\label{eq:chi} 
\chi^2(f_{Y,0}; f_{Y, 1}) \leq c_{14} h^{2m+2r-1}N^{-2m-1}. 
\end{equation}
\par 
The same upper bound holds for non--integer  $r$. 
Indeed, it follows from the above bounds on $I_1$ and $I_2$ that for integer $k\geq 1$  
\[
 \int_{-\infty}^\infty |(g\ast \psi_h)(x)|^2\rd x \leq \frac{c_{15} h^{2m+1}}{N^{2m-1}},\;\;
 \int_{-\infty}^\infty (1+ x^2)^k |(g\ast \psi_h)(x)|^2\rd x \leq \frac{c_{16} h^{2m-2k+1}}{N^{2m-1}},\;\;
\]
Then for real $0\leq r\leq k$ by the  interpolation inequality for the Sobolev spaces [see, e.g., \cite[Proposition~6.3.3]{Aubin00}]
we have 
\[
 \int_{-\infty}^\infty  (1+x^2)^r|(g\ast \psi_h)(x)|^2\rd x \leq \Big(\frac{c_{15} h^{2m+1}}{N^{2m-1}}\Big)^{1-r/k}
 \Big(\frac{c_{16} h^{2m-2k+1}}{N^{2m-1}}\Big)^{r/k} \leq c_{17} \frac{h^{2m-2r+1}}{N^{2m+1}},
\]
which 
yields (\ref{eq:chi}) for real $r$.
The choice 
\[
N=N_*:=\bigg(\frac{A}{h^{2r}}\bigg)^{1/\alpha},\;\;\;h=h_*:=\bigg(\frac{A^{(2m+1)/\alpha}}{n}\bigg)^{\frac{1}{2m+2r-1+2r(2m+1)/\alpha}}
\]
ensures that $f_0$ and $f_1$ are not distinguishable from observations which
leads to the lower bound
\begin{equation}\label{eq:rate-1}
 \cR_{n, \Delta_{x_0}}^*[\sF_{\alpha, p}(A, B)] \geq C \bigg(\frac{A^{(2m+1)/\alpha}}{n}\bigg)^{\frac{2r}{2m+2r-1+2r(2m+1)/\alpha}},\;\;p<2r-1. 
\end{equation}
 \par
 4$^0$. The rate of convergence obtained in  (\ref{eq:rate-1}) dominates the standard rate if 
 \[
  \frac{2r}{2r+2m-1+2r(2m+1)/\alpha}\leq \frac{\alpha}{2\alpha+2m+1} \;\;\Leftrightarrow\;\; 2r\leq 2m-1.
 \]
Therefore if 
$p<2r-1 \leq 2m-2$
then 
the standard rate of convergence is not achievable.
This completes the proof of (\ref{eq:low-b-1}).
\paragraph{Proof of statement (b)} 
Let $r>1/2$, and $f_0$ be defined in (\ref{eq:f-0}).
Clearly,  $f_0$ given by (\ref{eq:f-0}) satisfies 
$f_0\in \sS_\alpha (A)\cap \sM_p(B)$ for  $p<2r-1$ and large enough $A$ and $B$.
\par 
Let $\widehat{\psi}_0$ be the function satisfying conditions 
(i)--(iii) in the proof of statement~(a).
For positive integer $k$ define 
\[
 \widehat{\psi}_{h, k} (\omega):= \widehat{\psi}_0 \Big(\frac{\omega-\pi k/\theta}{h}\Big)
 +\widehat{\psi}_0\Big(\frac{\omega+\pi k/\theta}{h}\Big).
\]
Note that $\widehat{\psi}_{h,k}$ is even,  
supported on $[-\pi k/\theta-h, -\pi k/\theta+h]\cup 
[\pi k/\theta-h, \pi k/\theta+h]$, and 
$\psi_{h, k}$ and $\psi_{h, k^\prime}$, $k\ne k^\prime$ have disjoint supports
because $h<\pi/\theta$. Moreover, 
\[
 \psi_{h,k}(x)=\frac{1}{2\pi}\int_{-\infty}^\infty
 \widehat{\psi}_{h,k}(\omega) e^{i\omega x}\rd \omega = 2h \psi_0(hx) \cos \Big(\frac{\pi k x}{\theta}\Big),\;x\in \bR.
\]
Let $N\geq 1$ be an integer number.
Define the following family of functions 
\begin{equation*}
 f_w(x):= f_0(x)+ c_0M \varphi_{h, w}(x),\;\;\;
 \varphi_{h, w}(x):=\sum_{k=N+1}^{2N} w_k\, \psi_{h, k}(x),\;\;\;w\in \{0,1\}^N,
\end{equation*}
where $M>0$ and $c_0>0$ are real numbers.
\par 
1$^0$. First we demonstrate that $f_w(x)$, $w\in \{0,1\}^N$ is a probability density
from the class $\sS_\alpha(A)\cap \sM_p(B)$,
provided that  constants $M$, $h$ and $N$ are chosen in an appropriate way. 
 We have 
$\int_{-\infty}^\infty \varphi_{h, w}(x)\rd x=0$ because $\widehat{\psi}_{h,k}(0)=0$ for 
all $k=N+1, \ldots, 2N$.
Moreover, similarly to the proof of statement~(a), 
since $\psi_0$ is rapidly decreasing,  
\[
 |\varphi_{h,w}(x)| =\bigg|\sum_{k=N+1}^{2N} w_k \psi_{h, k}(x)\bigg|
 \leq c_1 |x|^{-2r} h^{-2r+1} N.
\]
Therefore, if we  let
\begin{equation}\label{eq:MM}
 M=h^{2r-1}N^{-1}
\end{equation}
then 
by an appropriate choice of  constant $c_0$ 
we can ensure that $c_0 M|\varphi_{h, w}(x)| \leq f_0(x)$ for all $x$. Thus, 
$f_w(x)$ is indeed a probability density for any $w\in \{0,1\}^N$. 
This also shows that $f_w\in \sM_p(B)$ with $p<2r-1$.
Furthermore,  
\begin{align*}
 M^2 &\int_{-\infty}^\infty
|\omega|^{2\alpha} \big|\widehat{\varphi}_{h, w}(\omega)\big|^2 \rd \omega 
\\
 & = M^2 \sum_{k=N+1}^{2N} 
\int_{-\infty}^\infty |\omega|^{2\alpha}
\bigg[ \Big|\widehat{\psi}_0\Big(\frac{\omega- \pi k/\theta}{h}\Big)\Big|^2 +
\Big|\widehat{\psi}_0\Big(\frac{\omega+ \pi k/\theta}{h}\Big)\Big|^2
\bigg] \rd \omega 
 \\
 & \leq c_2 M^2 h \sum_{k=N+1}^{2N} \int_{-\infty}^\infty
\Big|\frac{\pi k}{\theta} +\xi h\Big|^{2\alpha} |\widehat{\psi}_0(\xi)|^2 \rd \xi
\leq c_3 M^2h N^{2\alpha+1}= c_3 h^{4r-1} N^{2\alpha-1},
 \end{align*}
where the last expression follows from (\ref{eq:MM}).
Therefore, if 
\begin{equation}\label{eq:smoothness}
 c_3 h^{4r-1} N^{2\alpha-1} \leq A^2
\end{equation}
then $f_w\in \sS_\alpha (A)$ for any $w\in \{0,1\}^N$. 
\par\medskip 
2$^0$.  
By
the Varshamov--Gilbert lemma [see, e.g. \cite[Lemma~2.9]{Tsybakov}] 
there exists a subset 
$W=\{w^{(0)}, \ldots, w^{(J)}\}$ of $\{0, 1\}^N$ 
such that $w^{(0)}=(0,\ldots, 0)$, $J\geq 2^{N/8}$, and 
any pair of vectors $w, w^\prime \in W$ 
are
distinct in at least $N/8$ entries (the Hamming distance between $w$ and $w^\prime$ is at least $N/8$). In what follows we consider the family 
functions $\{f_w, w\in W\}$.
Clearly, for any $w, w^\prime \in W$ we have 
\begin{equation}\label{eq:rho}
\rho^2:= \big\|f_w - f_{w^\prime}\big\|^2_2 \geq c_4 M^2 N h = c_4 h^{4r-1} N^{-1}.
\end{equation}
\par\medskip 
3$^0$. Next, we bound the $\chi^2$--divergence  between distributions of observations
corresponding to $f_0$ and $f_{w}$, 
$w\in W\backslash w^{(0)}$.
As in the proof of statement~(a) we have 
\begin{align*}
& \chi^2(f_{Y,0}; f_{Y, w}) = c_0^2M^2
 \int_{-\infty}^\infty \frac{|(g\ast \varphi_{h,w})(x)|^2}{(g\ast f_0)(x)} \rd x
 \\
 & \leq c_5 M^2\int_{-\infty}^\infty  |(g\ast \varphi_{h,w})(x)|^2 \rd x + 
 c_5 M^2 \int_{-\infty}^\infty   x^{2r}|(g\ast \varphi_{h,w})(x)|^2 \rd x
 =: c_5I_1 + c_6I_2.
\end{align*}
Furthermore, 
\begin{eqnarray*}
&& I_1= \frac{M^2}{2\pi} \int_{-\infty}^\infty
 |\widehat{g}(\omega)|^2 |\widehat{\varphi}_{h, w}(\omega)|^2 \rd \omega 
 \\
 && \;\;\;=\frac{M^2}{\pi} 
 \sum_{k=N}^{2N} w_k \int_{-\infty}^\infty \bigg|\frac{\sin \theta\omega}{\theta \omega}
 \bigg|^{2m} \bigg|\widehat{\psi}_0\Big(\frac{\omega - \pi k/\theta}{h}\Big)\bigg|^2
 \rd \omega
 \\
 && \;\;\;\leq c_7 M^2 h^{2m+1} N^{-2m+1} = c_7 h^{4r+2m-1} N^{-2m-1}.
\end{eqnarray*}
An upper bound on $I_2$ is derived as in the proof of statement~(a). In particular, 
taking into account that $\widehat{\varphi}_{h, w}(\omega)$ is a sum of functions with disjoint
supports, and repeating the steps from (\ref{eq:I22}) to  (\ref{eq:I222}) 
we obtain 
\begin{eqnarray*}
 I_2 &\leq& c_8 \sum_{j=0}^r \frac{M^2}{h^{2r-2j}} \sum_{l=1}^j 
 \sum_{k=N}^{2N} 
 \int_{-\infty}^\infty \frac{|\sin \theta\omega|^{2m-2l}}{|\theta\omega|^{2m}}
 \bigg|\widehat{\psi}_0^{(r-j)}\Big(\frac{\omega-\pi k/\theta}{h}\Big)\bigg|^2\rd\omega
 \\
 &\leq & c_9 M^2 h^{2m-2r+1} N^{-2m+1} = c_9 h^{2r+2m-1} N^{-2m-1}.
\end{eqnarray*}
Combining these bounds on $I_1$ and $I_2$ we obtain 
\[
 \chi^2(f_{Y, 0}; f_{Y, w}) \leq c_{10} h^{2r+2m-1} N^{-2m-1}.
\]
\par 
4$^0$. To complete the proof we use 
Theorem~2.7 from \cite{Tsybakov}; see also Lemma~4 from \cite{GL14}.
In particular, this result implies that if $N$ and $h$ satisfy 
\begin{equation}\label{eq:chi-restriction}
 n h^{2r+2m-1} N^{-2m-1}=c_{11} N
\end{equation}
then for large $n$ one has 
$\cR^*_{n,\Delta_2}[\sG_{\alpha, p}(A, B)]\geq c_{12}\rho$, 
where $p<2r-1$, and $\rho$ is defined in (\ref{eq:rho}).
Now we choose $h$ and $N$ so that (\ref{eq:chi-restriction}) and 
(\ref{eq:smoothness}) are satisfied. To this end define
\[
 \kappa:= 2m+2r-1+ 
 \frac{2m+2}{2\alpha-1}(4r-1),\;\;\;\mu:=\frac{2\alpha}{2\alpha-1} (4r-1),
\]
and set 
\[
 N=\bigg(\frac{A^2}{h_*^{4r-1}}\bigg)^{1/(2\alpha-1)},\;\;\;
 h=h_*:= 
 \bigg(
 \frac{A^{2(2m+2)/(2\alpha-1)} }{n}\bigg)^{1/\kappa}, 
\]
With this choice (\ref{eq:smoothness}) and (\ref{eq:chi-restriction}) hold, and 
\begin{equation}\label{eq;rho}
 \rho^2 \geq c_4 h_*^{4r-1}N^{-1}= 
 c_4 A^{-2/(2\alpha-1)} h_*^{2\alpha (4r-1)/(2\alpha-1)}=c_5
 A^{\beta} 
 n^{-\mu/\kappa}, 
\end{equation}
where 
$\beta:= 2\big[(2m+2)(\mu/\kappa) -1\big]/(2\alpha-1)$. The rate of convergence obtained in
(\ref{eq:rho}) dominates the standard rate of convergence if 
\[
 \frac{\mu}{\kappa}= \frac{2\alpha(4r-1)}{(2\alpha-1)(2m+2r-1)+
 (4r-1)(2m+2)} \;\leq\; \frac{2\alpha}{2\alpha+2m+1},
\]
which is equivalent to $2r\leq 2m$. Therefore if 
$p<2r-1 \leq 2m-1$
then the standard rate of convergence is not attained. 
This completes the proof.
\epr

\paragraph{Acknowledgement}
The authors are grateful to Taeho Kim for careful reading and useful remarks. 
This article was prepared within the framework of the HSE University Basic Research
Program and funded by the Russian Academic Excellence Project ’5-100’.   AG is supported
by the ISF Reserarch Grant no.~361/15.

\bibliographystyle{imsart-nameyear}

\end{document}